\documentclass{amsart}

\newtheorem{theorem}{Theorem}[section]
\newtheorem{corollary}[theorem]{Corollary} 
\newtheorem{lemma}[theorem]{Lemma}
\newtheorem{example}[theorem]{Example}
\numberwithin{equation}{section}

\newcommand{\tref}[1]{Theorem~\textup{\ref{#1}}}
\newcommand{\cref}[1]{Corollary~\textup{\ref{#1}}}
\newcommand{\lref}[1]{Lemma~\textup{\ref{#1}}}

\newcommand{\C}{\mathbb C}

\def\new{\goodbreak\bigskip\noindent}

\def\J#1[#2]{\mathcal{J}_{#1\,}[#2]}
\def\JJ#1{\mathcal{J}_{#1}}
\def\B{\mathcal B}

\newcommand{\pol}{\mathcal P}
\newcommand{\con}{\,\mathcal C}

\newcommand{\upitem}[1]{\textup{(#1)}}

\newcommand{\xx}{x_1,x_2,x_3}
\newcommand{\yy}{y_1,y_2,y_3 }
\newcommand{\zz}{z_1,z_2,z_3 }
\newcommand{\uu}{u_1,u_2,u_3}
\newcommand{\vv}{v_1,v_2,v_3}

\begin{document}

\title{Completely Reducible Ternary Cubic Forms}
\author{Gary Brookfield}
\email{gbrookf@calstatela.edu}
\address{California State University, Los Angeles}
\begin{abstract}
We discuss various necessary and sufficient conditions for the complete reducibility of a ternary cubic form. In doing so, we prove the claim made in the 19th century that such a form is completely reducible if and only if its Hessian is a multiple of itself.
\end{abstract}

\maketitle

\section{Introduction}
The main goal of this article is to prove that a ternary cubic form \( f_{xxx} \) is completely reducible if and only if its Hessian \( \Delta_{xxx} \) is a multiple of \( f_{xxx} \) (\tref{thmq111}). See Sections~\ref{notation} and~\ref{cubic} for the notation and definitions. 

In one direction the claim is  easy to prove. If \( f_{xxx} \)  is completely reducible, that is, \( f_{xxx}=a_xb_xc_x \) for some linear forms \( a_x \), \( b_x \) and \( c_x \), then a direct calculation  gives \( \Delta_{xxx}= \J[a_x,b_x,c_x]^2\,a_xb_xc_x \) where \( \J[a_x,b_x,c_x] \) is the Jacobian of the linear forms, and so \( \Delta_{xxx} \) is a multiple of  \( f_{xxx} \). It turns out that the converse is quite hard to prove and has a curious history. 

Hessians of cubic forms were first investigated in 1844 by Otto Hesse \cite{hesse1844a,hesse1844b}. 
%
%
In 1846, Cayley \cite{cayley1846} was the first to notice that, if \( f_{xxx} \) is completely reducible, then its Hessian \( \Delta_{xxx} \) is a multiple of \( f_{xxx} \).
 In 1849, Aronhold \cite{aronhold1849}  stated explicitly, without proof,  that, if  \( \Delta_{xxx} \) is a multiple of \( f_{xxx} \), then \( f_{xxx} \) is completely reducible. (He expressed this result in terms of the \( 45 \)  \( 2\times2 \) minors of the \( 2\times10 \) matrix whose top row contains the coefficients of \( f_{xxx} \) and whose bottom row contains the coefficients of \( \Delta_{xxx} \).)
Subsequently, several mathematicians repeated the claim without proof, for example,  \cite[p.~182]{salmon1852} and \cite[p.~187]{sylvester1852}.  

In 1871, Gundelfinger \cite {gundelfinger1871b} collected many of the new discoveries about ternary cubic forms and attempted a proof that, if  \( \Delta_{xxx} \) is a multiple of \( f_{xxx} \), then \( f_{xxx} \) is completely reducible. Unfortunately, he used a dubious geometric argument to provide conditions for the reducibility of the form and ignored the case when \( \Delta_{xxx}=0 \). In a later paper from the same year, Gundelfinger \cite{gundelfinger1871c}  conceded the weakness of this argument and presented a longer proof that avoided the problematic geometric argument, but still failed to cover the case when \( \Delta_{xxx}=0 \). This same incomplete argument appears in Clebsch's comprehensive \textit{Vorlesungen \"uber Geometrie}  \cite[p.~597]{clebsch} in 1876. 

In the following decades, other mathematicians found conditions for the complete reducibility of ternary cubic forms without making any reference to the properties of the Hessian \cite{brill1898,brioschi1876,copeland,glenn}.  After about 1900, the claim that \( f_{xxx} \) is completely reducible if and only if \( \Delta_{xxx} \) is a multiple of \( f_{xxx} \) seems to have disappeared from the literature.

\new 
We will use 19th century  tools, such as transvectants and polars, that are no longer well-known, and so Sections~\ref{notation}-\ref{basis}  provide a condensed presentation of these concepts, without a discussion of their importance to the theory of invariants. Indeed, the invariant properties of forms are not discussed at all.

Our discussion depends on the validity of some hundred identities that relate the concomitants of forms. These are given without proof.  For example, the main theorem of the article uses the identity \eqref{eqq127}
\[ L_{xxx}+\delta f_{xxx}=[xyz]^2 \Vert \Gamma_{4ux} \Vert_{u\mapsto[yz]} \]
No proof of this identity is provided since all one has to do is use the definitions to expand both sides of the equation and confirm that they are equal. Because the expressions being equated have 85,119 terms when fully expanded, this is a job for a computer algebra system. 

This  article can be seen as an extensive elaboration of  \textit{Factoring Forms} \cite{brookfield} from 2016.

\section{Notation}\label{notation}

All forms in this article have coefficients in \( \C \), the field of complex numbers,
and almost all are ternary, that is, they are homogeneous polynomials in three variables \( \xx \). Frequently we will use other sets of covariant variables, \( \yy \) and \( \zz \). We will also use sets of contravariant variables, \( \uu \) and  \( \vv \). 
A geometer would view \((\xx) \), \( (\yy) \) and \( (\zz) \) as the coordinates of a point, and \( (\uu) \) and  \( (\vv) \) as the coordinates of a line. The point \( (\xx)\) is on the line \((\uu)\), equivalently, the line \( (\uu)\) passes through the point \((\xx)\),  if and only if \(u_1x_1+u_2x_2+u_3x_3=0\). The linear form on the left side of this equation is written \(u_x\) in accordance with much of 19th century literature and the notation used in this article. The words ``covariant'' and ``contravariant'' come from the theory of invariants, a dominating theme of nineteenth century algebra \cite{grace}.

\new We will use the notation \( f_x \), \( f_{xx} \) and \( f_{xxx} \) for ternary forms of degrees \( 1 \), \( 2 \), and \( 3 \) in the variables \( x_1\), \(x_2\), \(x_3 \) whose expansions are written as
\begin{equation}\label{eq00}
\begin{aligned}
f_x&=f_1x_1+f_2x_2+f_3x_3\\[3pt]
f_{xx}&=f_{11} x_{1}^2 + f_{12} x_{1} x_{2} + f_{22} x_{2}^2 + f_{13} x_{1} x_{3} + f_{23} x_{2} x_{3} +  f_{33}  x_{3}^2\\[3pt]
f_{xxx}&=  f_{111} x_{1}^3 + f_{112} x_{1}^2 x_{2} + f_{122} x_{1} x_{2}^2 + f_{222} x_{2}^3 +  f_{113} x_{1}^2 x_{3} + f_{123} x_{1} x_{2} x_{3} \\
&\qquad+ f_{223} x_{2}^2 x_{3} + f_{133} x_{1} x_{3}^2 +  f_{233} x_{2} x_{3}^2 + f_{333} x_{3}^3
\end{aligned}
\end{equation}
The coefficients are labeled with integer subscripts in increasing order. 


The symbol \( f \) without subscripts could represent a form of any degree. Partial derivatives with respect to \( \xx \) are abbreviated as follows: \[ \partial_1 f=\dfrac{\partial f}{\partial x_1}  \qquad \partial_2 f=\dfrac{\partial f}{\partial x_2} \qquad \partial_{12} f=\dfrac{\partial^2 f}{\partial x_1 \partial x_2} \qquad  \partial_{22} f=\dfrac{\partial^2 f}{\partial x_2\partial x_2}, \text{ etc.}\] 
For example, for any form \( f \) in \( \xx \), we have Euler's identity: 
\begin{equation}\label{eqw77}
(\deg f)\,f= x_1\partial_1f+x_2\partial_2f+x_3\partial_3f
\end{equation}

Similarly, 
\begin{align*}
F_u&=F_1u_1+F_2u_2+F_3u_3\\[3pt]
F_{uu}&=F_{11} u_{1}^2 + F_{12} u_{1} u_{2} + F_{22} u_{2}^2 + F_{13} u_{1} u_{3} + F_{23} u_{2} u_{3} +  F_{33}  u_{3}^2
\end{align*}
are forms in the variables \( \uu \), and \[ G_{ux}=G_{11}u_1x_1+ G_{12}u_1x_2+G_{13}u_1x_3+\cdots
+G_{32}u_3x_2+G_{33}u_3x_3\]
is a form in \( \xx \) and \( \uu \). 

For compactness we write
\begin{equation}\label{eqw61}
[xyz]=
\begin{vmatrix}
x_1 & x_2 & x_3\\
y_1 & y_2 & y_3\\
z_1 & z_2 & z_3 
\end{vmatrix} \qquad 
[abc]=
\begin{vmatrix}
a_1 & a_2 & a_3\\
b_1 & b_2 & b_3\\
c_1 & c_2 & c_3 
\end{vmatrix}
\end{equation}
where \( a_{x}=a_1x_1+a_2x_2+a_3x_3 \),  \( b_{x}=b_1x_1+b_2x_2+b_3x_3 \) and \( c_{x}=c_1x_1+c_2x_2+c_3x_3 \) are linear forms. 

 The most important tools in invariant theory are substitutions. If \( \sigma \) is a substitution and \( f \) is some form, then we write \( \Vert f\Vert_{\sigma} \) for \( f \) after the substitution \( \sigma \) has been made. For compactness, \( \sigma \) is written in symbolic form:
\begin{equation}\label{equ92}
\def\strut{\vrule height 10pt width 0pt depth 5pt}
\begin{array}{c|lcccccccc}
 \text{symbol}  & \text{substitution}\\\hline
\strut  x\mapsto y & x_1\mapsto y_1 \quad x_2\mapsto y_2 \quad x_3\mapsto y_3\\\hline
\strut u\mapsto [yz] &   u_1 \mapsto  y_2 z_3-y_3 z_2 \quad u_2\mapsto y_3 z_1 - y_1 z_3\quad  u_3 \mapsto  y_1 z_2-y_2 z_1\\\hline
\strut x\mapsto [uv] &   x_1 \mapsto  u_2 v_3-u_3 v_2 \quad x_2\mapsto u_3 v_1 - u_1 v_3\quad  x_3 \mapsto  u_1 v_2-u_2 v_1\\\hline
 \end{array}
\end{equation}
 For example, \(  \Vert u_x\Vert_{ u\mapsto [yz]}=[xyz] \). There are also substitutions that act on the coefficients of forms. If \( f \) and \( g \) are forms of the same degree, then the substitution \( f\mapsto g \) replaces the coefficients of \( f \) by the coefficients of \( g \). If \( f_{xx} \) and \( g_{xx}\) are quadratic forms, then \( \Vert f_{12}^2-4 f_{11}f_{22} \Vert_{f\mapsto g}=g_{12}^2-4 g_{11}g_{22}  \). 

Most important are substitutions that involve the derivative of forms. Let \( f \) be an arbitrary form. Then, for a linear form  \( a_x=a_1x_1+a_2x_2+a_3x_3 \), we define 
\begin{equation}\label{eq01}
\def\strut{\vrule height 10pt width 0pt depth 5pt}
 \begin{array}{l|l}
  \text{symbol}  & \text{substitution}\\\hline
 \strut a\mapsto \partial f & a_1   \mapsto \partial_1 f   \quad  a_2   \mapsto \partial_2 f   \quad   a_3   \mapsto \partial_3 f\\
 \hline
  a^2 \mapsto \partial^2 f &  
  \begin{aligned}
 \strut a_1^2    &\mapsto   \partial_{11}f   &   a_1a_2    &\mapsto  \partial_{12}f   &   a_1a_3    &\mapsto  \partial_{13}f   \\
 \strut a_2^2    &\mapsto  \partial_{22}f    &  a_2a_3   & \mapsto  \partial_{23}f    &  a_3^2    &\mapsto  \partial_{33}f
  \end{aligned} \\
 \hline
 \strut  a^3 \mapsto \partial^3 f &  
 \begin{aligned}
\strut  a_1^3 &\mapsto \partial_{111}f & a_1^2 a_2 &\mapsto \partial_{112}f & a_1 a_2^2 &\mapsto \partial_{122}f \\ a_2^3 &\mapsto \partial_{222}f &
 a_1^2 a_3 &\mapsto \partial_{113}f & a_1 a_2 a_3 &\mapsto \partial_{123}f \\
 a_2^2 a_3 &\mapsto \partial_{223}f & a_1 a_3^2 &\mapsto \partial_{133}f &
 a_2 a_3^2 &\mapsto \partial_{233}f \\ a_3^3 &\mapsto \partial_{333}f
\end{aligned}
\\\hline
 \end{array}
\end{equation}
%
%
%
%
For example, if a form \( f \) has degree \( n =\deg f\), then \[ \Vert a_x\Vert_{a\mapsto \partial f}=nf \qquad \Vert a_x^2\Vert_{a^2\mapsto \partial^2 f}=  n(n -1) f\qquad \Vert a_x^3\Vert_{a^3\mapsto \partial^3 f}=  n(n-1)(n-2) f.
 \] 
 The first equation is just a version of Euler's identity \eqref{eqw77}.

\section{Polar and Contraction Operators}

Introducing new variables \( X \), \( Y \) and \( Z \),  let \( x\mapsto XYZ \) denote  the substitution
\begin{equation}\label{eq33}
\begin{aligned}
x_1&\mapsto X x_1+Y y_1+Z z_1\\ 
x_2&\mapsto X x_2+Y y_2+Z z_2\\ 
x_3&\mapsto X x_3+Y y_3+Z z_3\\ 
\end{aligned}
\end{equation}
For a cubic form  \( f_{xxx} \),   \( \Vert f_{xxx}\Vert_{x\mapsto XYZ}  \)  is  a cubic form in \( X \), \( Y \) and \( Z \):
\begin{equation}\label{eqs11}
\begin{aligned}
\Vert f_{xxx}\Vert_{x\mapsto XYZ} &=f_{xxx}X^3+ f_{yyy}Y^3+f_{zzz}Z^3+ f_{xyz}XYZ+ f_{xxy}X^2Y\\
&\quad
+ f_{yyz} Y^2 Z+f_{xzz} X Z^2+f_{xyy}XY^2+f_{yzz} YZ^2+f_{xxz}X^2 Z
\end{aligned}
\end{equation}
The coefficients \( f_{xxy} \), \( f_{xyy} \) \( f_{xxy} \), \dots, \( f_{xzz} \) are called \textbf{polars} of \( f_{xxx} \). It is easy to see that \( f_{yyy}=\Vert f_{xxx}\Vert_{x\mapsto y} \) and \( f_{zzz}=\Vert f_{xxx}\Vert_{x\mapsto z} \). Other polars are more complicated, for example,
\begin{equation}\label{eqt66}
\begin{aligned}
f_{xxz} &= z_1 \partial_1f_{xxx}+z_2 \partial_2f_{xxx}+z_3 \partial_3f_{xxx}\\[5pt]
f_{xyz}&=  y_1 z_1 \,\partial_{11}f_{xxx}+  y_2 z_2 \,\partial_{22}f_{xxx}+ y_3 z_3  \,\partial_{33}f_{xxx}  + 
 (y_1 z_2 + y_2 z_1) \,\partial_{12}f_{xxx}  \\
&\qquad+  (y_1 z_3 + y_3 z_1) \,\partial_{13}f_{xxx}+ 
 (y_2 z_3 + y_3 z_2) \,\partial_{23}f_{xxx}
\end{aligned}
\end{equation}
The subscripts on  the forms \( f_{xxz} \) and \( f_{xyz} \) indicate the variables they contain and their degree. We will write the subscripts in alphabetic order.

Similarly, the polars of a linear form \( f_x \) and quadratic form \( f_{xx} \) are the coefficients of \( \Vert f_{x}\Vert_{x\mapsto XYZ} \) and \( \Vert f_{xx}\Vert_{x\mapsto XYZ} \):
\begin{align*}
 \Vert f_{x}\Vert_{x\mapsto XYZ}&= f_x X+ f_y Y+ f_z Z\\
\Vert f_{xx}\Vert_{x\mapsto XYZ} &=f_{xx} X^2+ f_{yy} Y^2+ f_{zz} Z^2+ f_{xy}XY+f_{yz} Y Z+ f_{xz} XZ
\end{align*}

The polars of a form are related to each other by \textbf{polar operators}, \( \pol^y_x \), \( \pol^z_x \),  \( \pol^x_y \),  \( \pol^y_z \), etc.  defined by 
\begin{equation}\label{eqg59}
 \pol_y^x[f]= y_1\frac{\partial f}{\partial x_1} +y_2\frac{\partial f}{\partial x_2}+y_3\frac{\partial f}{\partial x_3}
\end{equation}

For example,
\begin{align*}
f_{xxy}&=\pol^x_y[f_{xxx}]  & 2 f_{xyy}&=\pol^x_y[f_{xxy}] & 3 f_{yyy}&=\pol^x_y[f_{xyy}] \\
f_{xyy}&=\pol^y_x[f_{yyy}]  & 2 f_{xxy}&=\pol^y_x[f_{xyy}] & 3 f_{xxx}&=\pol^y_x[f_{xxy}]\\
f_{xyz}&=\pol^y_z[f_{xyy}]  & 2 f_{xzz}&=\pol^y_z[f_{xyz}] & 2 f_{xzz}&=\pol^y_x[f_{xyz}]
\end{align*}

The most important consequence of all these equations is that if any one of the polars is zero, then all the polars are zero. For example, from  \( 6f_{xxx}=\pol^y_x[\pol^z_x[f_{xyz}]] \) and \( f_{xyz}=\pol^x_y[\pol^x_z[f_{xxx}]]  \), we see that \( f_{xxx}=0 \) if and only if \( f_{xyz}=0 \).


\new Of course, the same ideas apply to forms in the contravariant variables: If \( F_{uuu} \) is a cubic form in \( \uu \), then \( F_{uuv}=\pol^u_v[F_{uuu}] \), \( F_{uvw}=\pol^u_v[\pol^u_w[F_{uuu}]]  \), etc. Also \( F_{uuu}=0 \) if and only if \( F_{uuv}=0 \), if and only if \( F_{uvw}=0 \), etc.

Important point: In this article, we will frequently, and without comment,  replace a hypothesis or conclusion, such as \( f_{xxx}=0 \) or \( F_{uuu}=0 \) or \( \theta_{uuxx}=0 \), by an equivalent condition, such as \( f_{xyz}=0 \) or \( F_{uvw}=0 \) or \( \theta_{uvzz}=0 \).



\new  The \textbf{contraction operator} \( \con_{ux}\) is defined by 
\[ \con_{ux}[f ]= \frac{\partial^2 f }{\partial x_1\partial u_1} +\frac{\partial^2 f }{\partial x_2\partial u_2}+\frac{\partial^2 f }{\partial x_3\partial u_3}. \]
for all  forms \( f \) in \( \xx \) and \( \uu \). 
For example, if \( F_{u}= F_{1} u_1 + F_{2} u_2 + F_{3} u_3\) and \( f_{xx} \) is as in \eqref{eq00}, then
\begin{multline*}
\con_{ux}[F_{u}f_{xx} ]= (2 F_1 f_{11} + f_{12} F_2+ f_{13} F_3) \,x_1 \\
+ (F_1 f_{12} + 2 F_2 f_{22} + f_{23} F_3) \,x_2 
+ (F_1 f_{13} + F_2 f_{23} + 2 F_3 f_{33}) \,x_3
\end{multline*}

The operators \( \con_{vx} \), \( \con_{uz} \), etc.~are defined similarly.
Each of these operators produces a new form whose degree in both the covariant and contravariant variables has been reduced by one. 


\section{Jacobians, Hessians and Transvectants}

The \textbf{Jacobian}  of ternary forms \( f \), \( g \) and \( h \) in \( \xx \) is  defined by
\[ \J[f,g,h]=
\begin{vmatrix}
\partial_1 f & \partial_2 f & \partial_3 f \\ \partial_1 g & \partial_2 g & \partial_3 g \\ \partial_1 h & \partial_2 h & \partial_3 h 
\end{vmatrix}\] 
For linear forms \( a_x \), \( b_x \) and \( c_x \), the Jacobian is simply the determinant of the coefficients:
\begin{equation}\label{eqq947}
\J[a_x,b_x,c_x]=
\begin{vmatrix}
a_1  &a_2  & a_3 \\ b_1  &b_2  & b_3 \\c_1  &c_2  & c_3 
\end{vmatrix}=[abc]
\end{equation}
Thus  \(\J[a_x,b_x,c_x]=0 \) if and only if \( \{a_x,b_x,c_x\} \) is linearly dependent. 

With this notation and \eqref{eq01}, the Jacobian of  \( f \), \( g \) and \( h \) can be expressed as 
\[  \J[f,g,h]=\Vert [abc]\Vert_{\substack{a\mapsto \partial f\\b\mapsto \partial g\\ c\mapsto \partial h}} \] Here we are extending our notation so that, if \( \sigma \) and \( \tau \) are substitutions that act on a form \( f \), then we write \( \Vert f\Vert_{\substack{\sigma\\\tau}}=\Vert \Vert f\Vert_\sigma\Vert_\tau =\Vert \Vert f\Vert_\tau\Vert_\sigma \) so long as the  second equality holds.

\new Other ways of combining three forms to create a new form are possible, in fact, the Jacobian is just the first \textbf{transvectant} of \( f \), \( g \) and \( h \). The second and third transvectants of \( f \), \( g \) and \( h \) are defined by
\[ \J2[f,g,h]=\Vert [abc]^2\Vert_{\substack{a^2\mapsto \partial ^2f\\b^2\mapsto \partial^2 g\\ c^2\mapsto \partial^2 h}} \qquad 
 \J3[f,g,h]=\Vert [abc]^3\Vert_{\substack{a^3\mapsto \partial ^3f\\b^3\mapsto \partial^3 g\\ c^3\mapsto \partial^3 h}} \] 
The higher order transvectants \( \J4[f,g,h] \), \( \J5[f,g,h] \), etc. can be defined similarly.

The calculations of these transvectants are straightforward. For example, since
\begin{align*}
 [abc]^2&=a_3^2 b_2^2 c_1^2 - 2 a_2 a_3 b_2 b_3 c_1^2 + a_2^2 b_3^2 c_1^2 - 2 a_3^2 b_1 b_2 c_1 c_2  + 2 a_2 a_3 b_1 b_3 c_1 c_2   \\
 &\qquad+ 2 a_1 a_3 b_2 b_3 c_1 c_2  -  2 a_1 a_2 b_3^2 c_1 c_2 + a_3^2 b_1^2 c_2^2 - 2 a_1 a_3 b_1 b_3 c_2^2  + a_1^2 b_3^2 c_2^2   \\
 &\qquad+ 2 a_2 a_3 b_1 b_2 c_1 c_3 - 2 a_1 a_3 b_2^2 c_1 c_3 -  2 a_2^2 b_1 b_3 c_1 c_3 + 2 a_1 a_2 b_2 b_3 c_1 c_3   \\
 &\qquad - 2 a_2 a_3 b_1^2 c_2 c_3 +  2 a_1 a_3 b_1 b_2 c_2 c_3 + 2 a_1 a_2 b_1 b_3 c_2 c_3 - 2 a_1^2 b_2 b_3 c_2 c_3  \\
 &\qquad + a_2^2 b_1^2 c_3^2 - 2 a_1 a_2 b_1 b_2 c_3^2 + a_1^2 b_2^2 c_3^2,
\end{align*} 
we have 
\begin{align*}
\J2[f,g,h] &=   \partial_{33} f\, \partial_{22} g\, \partial_{11} h  - 2 \partial_{23} f\, \partial_{23} g\, \partial_{11} h  +    \partial_{22} f\, \partial_{33} g\, \partial_{11} h   \\
&\quad  - 2 \partial_{33} f\, \partial_{12} g\, \partial_{12} h+    2 \partial_{23} f\, \partial_{13} g\, \partial_{12} h  + 2 \partial_{13} f\, \partial_{23} g\, \partial_{12} h   \\
&\quad -    2 \partial_{12} f\, \partial_{33} g\, \partial_{12} h  + 2 \partial_{23} f\, \partial_{12} g\, \partial_{13} h   -  2 \partial_{22} f\, \partial_{13} g\, \partial_{13} h    \\
&\quad - 2 \partial_{13} f\, \partial_{22} g\, \partial_{13} h  + 2 \partial_{12} f\, \partial_{23} g\, \partial_{13} h  + \partial_{33} f\, \partial_{11} g\, \partial_{22} h  \\
&\quad    -  2 \partial_{13} f\, \partial_{13} g\, \partial_{22} h  +  \partial_{11} f\, \partial_{33} g\, \partial_{22} h  - 2 \partial_{23} f\, \partial_{11} g\, \partial_{23} h  \\
&\quad    + 2 \partial_{13} f\, \partial_{12} g\, \partial_{23} h  +    2 \partial_{12} f\, \partial_{13} g\, \partial_{23} h  - 2 \partial_{11} f\, \partial_{23} g\, \partial_{23} h  \\
&\quad    +    \partial_{22} f\, \partial_{11} g\, \partial_{33} h  - 2 \partial_{12} f\, \partial_{12} g\, \partial_{33} h  +    \partial_{11} f\, \partial_{22} g\, \partial_{33} h
\end{align*}
for all forms \( f \), \( g \) and \( h \). The expression for \( \J3[f,g,h] \) has \( 54 \) terms --- too many to write out here. 

Because transvectants are defined using  derivatives of forms, they are multilinear operators. The \( n^{\text {th}} \) transvectant involves \( n^{\text{th}} \) order partial derivatives, so \( \J{n}[f,g,h] \) is zero  if any of \( \deg f \), \( \deg g \) and \( \deg h \) is less  than \( n \).  When nonzero, \( \J{n}[f,g,h] \) is a form of degree \( \deg f+\deg g+\deg h-3n \). 

If \( n \) is even, then \( \J{n}[f,g,h] \) is unchanged by any permutation of the forms \( f \), \( g \) and \( h \). If \( n \) is odd, then  \( \J{n}[f,g,h] \) is unchanged by any even permutation of the forms \( f \), \( g \) and \( h \), but changes sign by any odd permutation.  In particular,   if \( n \) is odd, then \( \J{n}[f,f,h]=0 \) and \( \J{n}[f,f,f]=0 \), and, more generally, \( \J{n}[f,g,h] =0\) if \( \{f,g,h\} \) is linearly dependent. 

If \( n \) is even, \( \J{n}[f,f,f] \) may not be zero. In fact, 
\begin{equation}\label{eq388}
\J2[f,f,f]=6\begin{vmatrix}
\partial_{11} f & \partial_{12} f & \partial_{13} f \\  \partial_{12} f & \partial_{22} f & \partial_{23} f \\ \partial_{13} f & \partial_{23} f & \partial_{33} f 
\end{vmatrix}
\end{equation}
The determinant of the second partial derivatives of \( f \) in this equation is called the \textbf{Hessian} of \( f \).
Another useful relationship between transvectants and determinants is
 \begin{equation}\label{eq382}
  \J2[f, f, u_x^2] =- 4
\begin{vmatrix}
\partial_{11}f & \partial_{12}f & \partial_{12}f & u_1 \\
 \partial_{12}f & \partial_{22}f & \partial_{23}f  & u_2 \\
  \partial_{13}f & \partial_{23}f & \partial_{33}f  & u_3\\
u_1 & u_2 & u_3 & 0
\end{vmatrix}
 \end{equation}
 
\new The following  identity makes it possible to ``factor out'' sufficiently high powers of linear forms from certain transvectants:
  \begin{equation}\label{eqa22}
 \J{k}[f,g\,a_x^m,a_x^n]=\binom{n}{k} \J{k}[f,g,a_x^k]\,a_x^{n+m-k}  \qquad \text{for }1\leq k\leq n
\end{equation}
The right side of this equation is zero if the degree of \( g \) is less than \( k \), so as special cases we get, for example, 
\begin{equation}\label{eqJ}
\begin{aligned}
\J[f,a_x^2,a_x]&=0& 
\J2[f,g_x a_x ,a_x^3]&=0 & 
\J2[f,g_x a_x ,a_x^2]&=0\\
\J[f,a_x^3 ,a_x]&=0 &
\J2[f,g_x a_x^2 ,a_x^2]&=0 & 
\J3[f, g_{xx}a_x ,a_x^3]&=0
\end{aligned}
\end{equation}
With a view to proof of \lref{lemu91} in the following section, we notice three easy consequences of these equations for a cubic form \( g_{xxx} \).  The first equation in the second row says that if  if  \( a_x^3 \) divides \( g_{xxx} \), then \( \J[f,g_{xxx},a_x]=0 \) for all forms \( f \), and in particular, \(  \J[g_{xxx},a_x,u_x]=0\).  The second equation in the same row says that, if  \( a_x^2 \) divides \( g_{xxx} \), then \( \J2[f,g_{xxx},a_x^2]=0 \) for all forms \( f \), and in particular, \(  \J2[g_{xxx},a_x^2,u_x^2]=0\). Similarly, the last equation implies that, if \( a_x \) divides \( g_{xxx} \), then \( \J3[f,g_{xxx},a_x^3]=0 \) for all forms \( f \), and in particular, \( \J3[g_{xxx},a_x^3,u_x^3]=0 \). The converses of these claims are proved in \lref{lemu91}.

There are other identities that don't fit the pattern of \eqref{eqa22}:
\begin{equation}\label{eqL}
\begin{gathered}
2 \J2[f, a_x b_x, a_x b_x] =-\J2[f, a_x^2, b_x^2]\\
\J3[f,a_x^2 b_x,a_x^2 c_x]=0 \qquad
\J3[f,a_xb_x^2,a_x c_x^2]=\J3[f,b_xc_x^2,b_x a_x^2]
\end{gathered}
\end{equation}
In addition, almost directly from the definition, we get
\begin{equation}\label{eqK}
 \J{k}[a_x^k,b_x^k,c_x^k]=(k!)^3 \J[a_x,b_x,c_x]^k= (k!)^3 [abc]^k
\end{equation}

%
%

%

\new The forms that can be constructed from \( f \) and \( u_x  \) (and \( v_x \), etc.) using transvectants, polar operators, contraction operators, substitutions, etc.~are \textbf{concomitants}  of \( f \). For example,
\[ 
 f_{xyz} \quad 
 \con_{ux}[f_{xxx} \J[f_{xxx},f_{xxx},u_x^2]] \quad
 \J3[f_{xxx} u_x, f_{xxx} u_x,f_{xxx} u_x]\] 
 are concomitants of \( f_{xxx} \). Concomitants that are forms only in the covariant variables, \( \xx \), \( \yy \), \dots  are called \textbf{covariants} of \( f \). Concomitants that are forms only in the contravariant variables \( \uu \), \( \vv \), \dots  are called \textbf{contravariants} of \( f \). Concomitants in which none of the variables appear are called \textbf{invariants} of \( f \). For example, \( \J2[f_{xx},f_{xx},f_{xx}] \) is an invariant of \( f_{xx} \). For the precise definition of these terms and for their importance to the theory of invariants, see Grace and Young, \textit{The Algebra of Invariants}  \cite{grace}.

The idea in this section of using substitutions to construct concomitants of forms dates back to Aronhold \cite{aronhold1858} in 1858. From it evolved the symbolic method for invariant theory used extensively by 19th century mathematicians \cite{clebsch,osgood}.

\section{Bases for Forms}\label{basis}

Suppose that \( a_x \), \( b_x \)  and \( c_x \) are linear forms such that \( \J[a_x,b_x,c_x]=[abc] \) is nonzero. Then \( \mathcal{B}_1=\{a_x,b_x,c_x\} \)  is linearly independent, so forms a basis for the vector space of linear forms. Not too surprisingly, 
\begin{gather*}
\B_2=\{a_x^2,b_x^2,c_x^2,a_x  b_x,a_x  c_x , b_xc_x \}\\
\B_3=\{a_x^3,b_x^3,c_x^3,a_x^2 b_x,a_x^2 c_x , b_x^2 a_x,b_x^2 c_x, c_x^2 a_x, c_x^2 b_x,a_x b_x c_x\}
\end{gather*}
are bases for the vector spaces of quadratic forms and of cubic forms. 

This is made clear  by the following identities that express arbitrary  forms \( f_x \),  \( f_{xx} \) and \( f_{xxx} \) as  explicit linear combinations of  forms in \( \B_1 \), \( \B_2 \) and \( \B_3 \). 

\begin{align}\label{equ13}
 [abc] f_x &= \J[f_x, b_x, c_x]\, a_x + \J[f_x, c_x, a_x]\, b_x + \J[f_x, a_x, b_x]\, c_x \\[7pt]
\begin{split}\label{equ23}
2^3 [abc]^2 f_{xx} &= \J2[f_{xx}, b_x^2, c_x^2] a_x^2 
+  \J2[f_{xx}, a_x^2, c_x^2] b_x^2 
+ \J2[f_{xx}, a_x^2, b_x^2] c_x^2\\
&\quad- 2 \J2[f_{xx}, a_x b_x, c_x^2] a_x b_x 
-2 \J2[f_{xx}, a_x c_x, b_x^2] a_x c_x\\
&\quad -2 \J2[f_{xx}, b_x c_x, a_x^2] b_x c_x
\end{split}\\[7pt]
\begin{split}\label{equ33}
 6^3 [abc]^3 f_{xxx} &=  
 \J3[f_{xxx}, b_x^3 , c_x^3]\, a_x^3 
+ \J3[f_{xxx}, c_x^3 , a_x^3]\, b_x^3 \\
&\quad+ \J3[f_{xxx}, a_x^3 , b_x^3]\, c_x^3 
+ 3 \J3[f_{xxx}, c_x^3, b_x^2 a_x]\, a_x^2 b_x \\
&\quad+ 3\J3[f_{xxx}, b_x^3, a_x^2 c_x]\, c_x^2 a_x 
 +3 \J3[f_{xxx}, a_x^3, c_x^2 b_x]\, b_x^2 c_x \\
&\quad - 3\J3[f_{xxx}, a_x^3, b_x^2 c_x]\, c_x^2 b_x 
 - 3\J3[f_{xxx}, b_x^3, c_x^2 a_x]\, a_x^2 c_x \\
& \quad - 3\J3[f_{xxx}, c_x^3, a_x^2 b_x]\, b_x^2 a_x
 + 9 \J3[f_{xxx}, a_x b_x^2, a_x c_x^2]\, a_x b_x c_x 
\end{split}
\end{align}
Because of \eqref{eqL} we have 
 \[ \J3[f_{xxx}, a_x b_x^2, a_x c_x^2]=\J3[f_{xxx}, b_x c_x^2, b_x a_x^2]=\J3[f_{xxx}, c_x a_x^2, c_x b_x^2], \]
and so the coefficient of \( a_x b_x c_x \)  in \eqref{equ33}  has the expected symmetry with respect to permutations of \( \{a_x,b_x,c_x\} \).


\new Certain concomitants can also be expressed as linear combinations of forms in  \( \B_1 \), \( \B_2 \) and \( \B_3 \). For example, 
%
%
\begin{multline} \label{eqv33}
4  [abc]\, \J[f_{xx}, a_x, b_x] =-   \J2[f_{xx}, a_x c_x, b_x^2]\, a_x\\
 -   \J2[f_{xx}, b_x c_x, a_x^2]\, b_x +  \J2[f_{xx}, a_x^2, b_x^2]\, c_x 
\end{multline}
\begin{multline}\label{eqv37} 
6^3 [abc]^2 \J[f_{xxx}, a_x, b_x] = 
3 \J3[f_{xxx}, a_x^3, b_x^3]\, c_x^2 + 6 \J3[f_{xxx}, b_x^3, a_x^2 c_x]\, a_x c_x \\
+ 3 \J3[f_{xxx}, a_x^3, c_x^2 b_x ]\, b_x^2 - 6 \J3[f_{xxx}, a_x^3, b_x^2 c_x]\, b_x c_x \\
- 3 \J3[f_{xxx}, b_x^3, c_x^2 a_x]\, a_x^2 + 9 \J3[f_{xxx}, a_x b_x^2, a_x c_x^2]\, a_x b_x
\end{multline}
\begin{multline} \label{eqv38}
9 [abc]\, \J2[f_{xxx}, a_x^2, b_x^2] =\J3[f_{xxx}, b_x^3, a_x^2 c_x]\, a_x  \\
- \J3[f_{xxx}, a_x^3, b_x^2 c_x]\,  b_x + \J3[f_{xxx}, a_x^3, b_x^3]\, c_x
\end{multline}
\begin{multline}\label{eqv39}
18 [abc]\, \J2[f_{xxx}, a_x^2, b_x c_x] =
-3 \J3[f_{xxx}, a_x b_x^2, a_x c_x^2]\, a_x \\
- 2  \J3[f_{xxx}, a_x^3, c_x^2 b_x] b_x + 2 \J3[f_{xxx}, a_x^3, b_x^2 c_x]\, c_x
\end{multline}

%

An immediate consequence of  \eqref{equ13}  is that, if \( \{a_x,b_x\} \) is linearly  independent, then \( \J[f_x,a_x,b_x]=0 \) if and only if \( f_x=a_0a_x+b_0b_x\) for some \( a_0,b_0\in \C \). There are similar results for quadratic and cubic forms:

\begin{lemma}\label{lemq913}
Let  \( a_x \) and \( b_x \) be  linear forms such that \( \{a_x,b_x\} \) is linearly independent.

\noindent Let  \( f_{x} \)  be a linear form.
\begin{enumerate}
\item   [\textup{(A1)}] \( \J[f_x,a_x,b_x]=0 \) if and only if \( f_x=g_0a_x+h_0b_x\) for some \( g_0,h_0\in \C \).
\end{enumerate}
Let  \( f_{xx} \)  be a quadratic form.
\begin{enumerate}
\item   [\textup{(B1)}]  \( \J[f_{xx},a_x, b_x] =0\) if and only if \( f_{xx} =g_0 a_x ^2+ h_0a_xb_x +   k_0 b_x^2 \) for some \( g_0,h_0,k_0\in \C \).
\item  [\textup{(B2)}]  \( \J2[f_{xx}, a_x^2, b_x^2] =0 \) if and only if \( f_{xx}= g_xa_x + h_xb_x \) for some linear forms \( g_x \) and \( h_x \).
\end{enumerate}
Let \( f_{xxx} \)  be a cubic form.
\begin{enumerate}
\item  [\textup{(C1)}] \( \J[f_{xxx},a_x, b_x] =0\) if and only if  \( f_{xxx} =g_0 a_x ^3+ h_0 a_x^2b_x+k_0a_xb_x^2+  i_0 b_x^3\) for some \( g_0,h_0,k_0,i_0\in \C \).

\item  [\textup{(C2)}] \( \J2[f_{xxx},a_x^2, b_x^2] =0\) if and only if  \( f_{xxx} =g_x a_x ^2+ h_x a_xb_x+k_xb_x^2\) for some linear forms \( g_x \), \( h_x \)  and \( k_x \).

\item  [\textup{(C3)}] \(\J3[f_{xxx}, a_x^3 , b_x^3] =0 \) if and only if \( f_{xxx}= g_{xx}a_x+ h_{xx} b_x\) for some quadratic forms \( g_{xx} \) and \( h_{xx} \).
\end{enumerate}

%
%
%
%
\end{lemma}

\begin{proof}
In all cases, if \( f_x \), \( f_{xx} \), \( f_{xxx} \) has the form on the right, then the transvectant on the left is zero because of the multilinearity of transvectants and  \eqref{eqJ}.

For the converses of these claims, choose a linear form \( c_x \) so that \(\B_1= \{a_x,b_x,c_x\} \) is linearly independent, that is,  \( [abc]\neq 0 \).

\noindent\upitem{A1} This follows directly from  \eqref{equ13}.

\noindent\upitem{B1}  If  \( \J[f_{xx},a_x, b_x] =0\), then, because of  \eqref{eqv33}  and the linear independence of \( \B_1 \), we get 
\(  \J2[f_{xx}, a_x c_x, b_x^2]=
  \J2[f_{xx}, b_x c_x, a_x^2]= \J2[f_{xx}, a_x^2, b_x^2] =0\).
Hence the coefficients of \( a_x c_x \), \( b_x c_x \) and  \( c_x^2 \) in \eqref{equ23} are zero and  \( f_{xx} \) has the claimed form.
 
\noindent\upitem{B2} If  \( \J2[f_{xx}, a_x^2, b_x^2] =0 \), then the coefficient of \( c_x^2 \) in \eqref{equ23} is zero. All other terms  of \eqref{equ23} are multiples of \( a_x \) or of \( b_x \)  and so  \( f_{xx} \) has the claimed form.
 
\noindent\upitem{C} The proof is similar to the proof of (B). In each case, the assumption about a transvectant of \( f_{xxx} \), \( a_x \) and \( b_x \) implies that certain coefficients on the right side of \eqref{equ33} are zero. In the cases (C1) and (C2), the particular coefficients are determined by the linear independence of \( \B_1 \) and \( \B_2 \),  and the identities \eqref{eqv38} and \eqref{eqv39}. Once those coefficients are set to zero, \eqref{equ33}  can be solved for \( f_{xxx} \), putting it into the claimed form.
\end{proof}

 We next investigate the conditions \( \J[f,a_x,u_x]=0 \), \( \J2[f,a_x^2,u_x^2]=0 \) and \( \J3[f,a_x^3,u_x^3]=0 \) where \( f \) is any form and \( a_x \) is a  linear form. First we notice that for any quadratic form \( g_{xx} \) we have  \( 2\J2[f,a_x^2,g_{xx}]=\Vert \J2[f,a_x^2,u_x^2] \Vert_{u^2\mapsto g} \), and so the hypothesis that \( \J2[f,a_x^2,u_x^2]=0 \) implies \( \J2[f,a_x^2,g_{xx}]=0 \) for all \( g_{xx} \) and, in particular, we can replace \( u_x^2 \) by any form in \( \B_2 \) to get \( \J2[f,a_x^2,b_x^2]= \J2[f,a_x^2,c_x^2]= \J2[f,a_x^2,b_x c_x]=0 \). Of course, a similar argument applies to the hypotheses \( \J[f,a_x,u_x]=0 \) and \( \J3[f,a_x^3,u_x^3]=0 \). 
 
%
%
%
%
%
 
\begin{lemma}\label{lemu91}
Let \( a_x \) be a nonzero linear form.

\noindent Let  \( f_x \) be a linear form.
\begin{enumerate}
\item [\textup{(A1)}] \( \J[f_{x},a_x,u_x]=0 \) if and only if \( a_x \) divides \( f_{x} \).
\end{enumerate}
Let  \( f_{xx} \) be a quadratic form.
\begin{enumerate}
\item [\textup{(B1)}] \( \J[f_{xx},a_x,u_x]=0 \) if and only if \( a_x^2 \) divides \( f_{xx} \).
\item [\textup{(B2)}] \( \J2[f_{xx},a_x^2,u_x^2]=0 \) if and only if \( a_x \) divides \( f_{xx} \).
\end{enumerate}
Let \( f_{xxx} \)  be a cubic form.
\begin{enumerate}
\item  [\textup{(C1)}]  \( \J[f_{xxx},a_x,u_x]=0 \) if and only if \( a_x^3 \) divides \( f_{xxx} \).
\item [\textup{(C2)}] \( \J2[f_{xxx},a_x^2,u_x^2]=0 \) if and only if \( a_x^2 \) divides \( f_{xxx} \).
\item [\textup{(C3)}] \( \J3[f_{xxx},a_x^3,u_x^3]=0 \) if and only if \( a_x \) divides \( f_{xxx} \).
\end{enumerate}

\end{lemma}

\begin{proof}
In all cases, if \( f_x \), \( f_{xx} \), \( f_{xxx} \) has property on the right, then the transvectant on the left is zero because of  \eqref{eqJ}.

To prove the converse, choose linear forms \( b_x \) and  \( c_x \) so that \(\B_1= \{a_x,b_x,c_x\} \) is linearly independent, that is,  \( [abc]\neq 0 \).

\noindent\upitem{A1} Suppose  that \( \J[f_{x},a_x,u_x]=0 \). As pointed out above, we can replace \( u_x\) by forms in \( \B_1 \) to get \( \J[f_{x},a_x,b_x]=0 \) and  \( \J[f_{x},a_x,c_x]=0 \). Hence  \eqref{equ13} becomes \( [abc] f_x= \J[f_x,b_x,c_x] \,a_x \),  showing that    \( a_x \) divides \( f_{x} \)


\noindent\upitem{B1} Suppose that \( \J[f_{xx},a_x,u_x]=0 \). Then, in particular, \( \J[f_{xx},a_x,b_x]=0 \) and so, by \eqref{eqv33} and the linear independence of \( \B_1 \), we get \( \J2[f_{xx},a_xc_x,b_x^2]=\J2[f_{xx},b_xc_x,a_x^2]=\J2[f_{xx},a_x^2,b_x^2]=0 \). Similarly, from  \( \J[f_{xx},a_x,c_x]=0 \) we get \( \J2[f_{xx},a_xb_x,c_x^2]=\J2[f_{xx},b_xc_x,a_x^2]=\J2[f_{xx},a_x^2,c_x^2]=0 \). Hence all terms on the right side of \eqref{equ23} that are not multiples of \( a_x^2 \) have zero coefficients. This  shows that   \( a_x^2 \) divides \( f_{xx} \).

\noindent\upitem{B2} Suppose that \( \J2[f_{xx},a_x^2,u_x^2]=0 \). Replacing  \( u_x^2\) by forms in \(  \B_2 \), we get \( \J2[f_{xx}, a_x^2, b_x^2] =\J2[f_{xx}, a_x^2, b_x c_x] =\J2[f_{xx}, a_x^2, c_x^2] = 0  \).  Hence all terms on the right side of \eqref{equ23} that are not multiples of \( a_x \) have zero coefficients. This  shows that  \( a_x \) divides \( f_{xx} \).

\noindent\upitem{C1} Suppose that \( \J[f_{xxx},a_x,u_x]=0 \). Replacing \( u_x \) by forms in \( \B_1 \) we get \( \J[f_{xxx},a_x,b_x]=\J[f_{xxx},a_x,c_x]=0  \).  Because  of  \( \J[f_{xxx},a_x,b_x]=0  \),  \eqref{eqv37} and the linear independence of \( \B_2 \), we see that six of terms on the right side of \eqref{equ33} are zero. Similarly  \( \J[f_{xxx},a_x,c_x]=0 \), and so three additional terms on the right side of \eqref{equ33} are zero. Hence only the first term can have a zero coefficient. This shows  that  \( a_x^3 \) divides \( f_{xxx} \).

\noindent\upitem{C2} Suppose that \( \J2[f_{xxx},a_x^2,u_x^2]=0 \).  
Replacing  \( u_x^2\) by forms in \( \B_2 \) we get  \( \J2[f_{xxx},a_x^2,b_x^2]=\J2[f_{xxx},a_x^2,b_xc_x]= \J2[f_{xxx},a_x^2,c_x^2]=0 \).  Because of \eqref{eqv38},  \eqref{eqv39} and the linear independence of \( \B_1 \), these equations imply that all terms on the right side of \eqref{equ33} that are not multiples of \( a_x^2 \) have zero coefficients.  This shows that \( a_x^2 \) divides \( f_{xxx} \).

\noindent\upitem{C3} Suppose that \( \J3[f_{xxx},a_x^3,u_x^3]=0 \). Replacing  \( u_x^3 \) by forms in \( \B_3 \), we see that the four terms on the right side of \eqref{equ33} that are not multiples of \( a_x \) have zero coefficients. This shows that \( a_x \) divides \( f_{xxx} \).
\end{proof}

Though hardly a surprise, (A1) of this lemma can be interpreted as saying that \( \{a_x, f_x\} \) is linearly dependent if and only if \( \J [f_x, a_x, u_x] = 0 \). Alternatively, replacing \( f_x \) by \( b_x \) and using \eqref{eqq947}, \( \{a_x, b_x\} \) is linearly dependent if and only if \(  [abu] = 0 \).

Another way of proving this lemma would be to use  identities  that relate transvectants of different orders:
 \begin{align}
\J2[f, a_x v_x, u_x^2] &=   4 \J[\J[f, a_x, u_x], v_x, u_x]\label{eqg67}\\
\J3[f, a_x v_x^2, u_x^3] &=   9 \J[\J2[f, a_x v_x, u_x^2], v_x, u_x] \notag \\
 \J3[f,a_x^2 v_x,u_x^3] &=9\J[\J2[f,a_x^2,u_x^2],v_x,u_x]] \notag
 \end{align}
If, for example, \( \J[f_{xx},a_x,u_x]=0 \), then from \eqref{eqg67}, we get \(  \J2[f_{xx}, a_x v_x, u_x^2]=0 \) for all \( v_x \) and \( u_x \). Replacing \( v_x \) by  forms in \( \B_1 \) and  \( u_x^2 \) by  forms in \( \B_2 \)  gives 
\begin{gather*}
\J2[f_{xx}, a_x^2, b_x^2]=  \J2[f_{xx}, a_x c_x, b_x^2]=\J2[f_{xx}, a_x^2, c_x^2]=0\\
 \J2[f_{xx}, a_x b_x, c_x^2]= \J2[f_{xx}, a_x c_x, b_x^2]=0 
\end{gather*}
So all terms on the right side of \eqref{equ23}, except the first,  are zero  and \( a_x^2 \) divides \( f_{xx} \), as claimed in (B1).

To demonstrate yet another way of proving \lref{lemu91},  we reprove   (C3).  Applying \( \con_{ux}^3 \) to \( \J3[f_{xxx},a_x^3,u_x^3]\), multiplied by  sufficient powers of \(b_x \) and \(c_x \), gives a form that can be written in the basis \( \B_3\):  
\begin{multline*}
\con_{ux}^3[\J3[f_{xxx}, a_x^3, u_x^3]\, b_x^3 c_x^3] =
 36 \Big( \J3[f_{xxx}, a_x^3, b_x^3]\,  c_x^3- \J3[f_{xxx}, c_x^3, a_x^3]\,  b_x^3 \\
 +   9 \J3[f_{xxx}, a_x^3, c_x^2 b_x] \, b_x^2 c_x 
 +   9 \J3[f_{xxx}, a_x^3, b_x^2 c_x]\,  c_x^2 b_x\Big)
\end{multline*}
If  \( \J3[f_{xxx},a_x^3,u_x^3]=0 \), then, because of the linear dependence of  \(\B_3\), we get
\begin{align*}
\J3[f_{xxx}, a_x^3, b_x^3] = \J3[f_{xxx}, c_x^3, a_x^3] &=0\\\J3[f_{xxx}, a_x^3, c_x^2 b_x] =\J3[f_{xxx}, a_x^3, b_x^2 c_x] &=0
\end{align*}   
  Thus the four terms on the right side of \eqref{equ33} that are not multiples of \( a_x \) have zero coefficients. This shows that \( a_x \) divides \( f_{xxx} \).
 
 We will use this method in all remaining proofs in this section.

\begin{lemma}\label{lemu46}
Let  \( a_x \) and \( b_x \) be  linear forms such that \( \{a_x,b_x\} \) is linearly independent. 

%
%
%
%
%

\begin{enumerate}
\item [\textup{(B)}] For a quadratic form \( f_{xx} \)  the following are equivalent:

\begin{enumerate}
\item [\textup{(1)}] \( \J2[f_{xx},a_x b_x,u_x^2] =0\)

\item [\textup{(2)}] \( f_{xx}=a_0 a_x^2+ b_0 b_x^2 \) for some \( a_0,b_0\in \C \).
\end{enumerate}

\item [\textup{(C)}] For a cubic form \( f_{xxx} \)  the following are equivalent:

\begin{enumerate}
 \item [\textup{(1)}]\( \J2[f_{xxx},a_xb_x,u_x^2]=0 \)
 
 \item [\textup{(2)}] \( f_{xxx}= a_0a_x^3+b_0 b_x^3\) for some \( a_0,b_0\in \C \).
 \end{enumerate}
\end{enumerate}

\end{lemma}

\begin{proof}
If \( f_{xx} \) or  \( f_{xxx} \) satisfies (2),  then (1) follows from \eqref{eqJ}.

For the converses of these claims, choose a linear form \( c_x \) so that \(\B_1= \{a_x,b_x,c_x\} \) is linearly independent, that is,  \( [abc]\neq 0 \).


\noindent\upitem{B} If \( \J2[f_{xx},a_x b_x,u_x^2] =0\), then the identity
\begin{multline*}
\con_{ux}^2[\J2[f_{xx}, a_x b_x, u_x^2] a_x b_x c_x^2] \!= 4 \J2[f_{xx}, a_x b_x, c_x^2] a_x b_x-2 \J2[f_{xx}, a_x^2, b_x^2] c_x^2 \\ 
 - 4  \J2[f_{xx}, a_x c_x, b_x^2] a_x c_x - 4  \J2[f_{xx}, b_x c_x, a_x^2] b_x c_x
\end{multline*}
and the linear independence of \( \mathcal{B}_2 \) imply that  the last  four terms on the right side of \eqref{equ23} are zero, that is,
\[ 8 [abc]^2 f_{xx} = \J2[f_{xx}, b_x^2, c_x^2]\, a_x^2 
+  \J2[f_{xx}, a_x^2, c_x^2] \,b_x^2 \] Since \(  [abc]\neq0 \), \( f_{xx} \) can be written in the claimed form.

\noindent\upitem{C} If  \( \J2[f_{xxx},a_xb_x,u_x^2]=0 \), then the identity
 \begin{align*}
 \con_{ux}^3[&\J[\J2[f_{xxx}, a_x b_x, u_x^2]  \,u_x c_x, a_x ^2 c_x , b_x^2 c_x]]\\
&= 24 \Big(-2 \J3[f_{xxx}, a_x^3, b_x^3]\, c_x^3 -  
    \J3[f_{xxx}, c_x^3, b_x^2 a_x]\, a_x^2 b_x  \\
&\qquad + 5 \J3[f_{xxx}, b_x^3, a_x^2 c_x]\, c_x^2 a_x + 
    \J3[f_{xxx}, a_x^3, c_x^2 b_x]\, b_x^2 c_x \\
&\qquad  -   5 \J3[f_{xxx}, a_x^3, b_x^2 c_x]\, c_x^2 b_x - 
    \J3[f_{xxx}, b_x^3, c_x^2 a_x]\, a_x^2 c_x  \\
&\qquad +  \J3[f_{xxx}, c_x^3, a_x^2 b_x]\, b_x^2 a_x - 
    12 \J3[f_{xxx}, a_x b_x^2, a_x c_x^2]\, a_x b_x c_x\Big)
\end{align*}
and the linear independence of \( \mathcal{B}_3 \) imply that the last eight terms of the right side of \eqref{equ33} are zero, that is,
\begin{equation}\label{eqq949}
216  [abc]^3 f_{xxx} =  
 \J3[f_{xxx}, b_x^3 , c_x^3]\, a_x^3 
+ \J3[f_{xxx}, c_x^3 , a_x^3]\, b_x^3. 
\end{equation}
Since \(  [abc]\neq0 \), \( f_{xxx} \) can be written in the claimed form. 
\end{proof}

\begin{lemma} 
 Let \( f_{xxx} \) be a cubic form, and let  \( a_x \), \( b_x \) and \( c_x \) be linear forms such that \( \{a_x,b_x,c_x\} \)  is linearly independent. The following are equivalent:
 \begin{enumerate}
 \item \( \J3[f_{xxx},a_xb_xc_x,u_x^3]=0 \)
 
 \item \( f_{xxx}=a_0 a_x^3+ b_0 b_x^3+c_0 c_x^3+d_0 a_x b_x c_x \) for some \( a_0,b_0,c_0,d_0\in \C \).
 \end{enumerate}
 \end{lemma}
 
 \begin{proof} 
 If  \( f_{xxx} \) satisfies (2),  then (1) follows from \eqref{eqJ} and the antisymmetry of \( \JJ3 \).
 
If \( \J3[f_{xxx},a_xb_xc_x,u_x^3]=0 \), then the identity
\begin{align*}
 \con_{ux}^3[&\J3[f_{xxx},  a_x b_x c_x, u_x^3]\, a_x^2 b_x^2 c_x^2] =\\
 & 24 \Big( \J3[f_{xxx}, c_x^3, b_x^2 a_x]\, a_x^2 b_x 
+ \J3[f_{xxx}, b_x^3, a_x^2 c_x]\, c_x^2 a_x \\
&\quad+ \J3[f_{xxx}, a_x^3, c_x^2 b_x]\, b_x^2 c_x 
+  \J3[f_{xxx}, a_x^3, b_x^2 c_x]\, c_x^2 b_x \\
 &\quad+  \J3[f_{xxx}, b_x^3, c_x^2 a_x]\, a_x^2 c_x + \J3[f_{xxx}, c_x^3, a_x^2 b_x]\, b_x^2 a_x\Big)
 \end{align*}
and the linear independence of \( \mathcal{B}_3 \) imply that six of the terms on the right side of \eqref{equ33} are zero. 
Then \eqref{equ33} becomes 
\begin{equation}
\begin{aligned}
216  [abc]^3 f_{xxx} &=  
 \J3[f_{xxx}, b_x^3 , c_x^3]\, a_x^3 
+ \J3[f_{xxx}, c_x^3 , a_x^3]\, b_x^3 \\
&\quad+ \J3[f_{xxx}, a_x^3 , b_x^3]\, c_x^3 
 + 9 \J3[f_{xxx}, a_x b_x^2, a_x c_x^2]\, a_x b_x c_x .
\end{aligned}
\end{equation}
Since \(  [abc]\neq0 \), \( f_{xxx} \) can be written in the claimed form. 
 \end{proof}

 \begin{lemma}
 Let \( f_{xxx} \) be a cubic form, and let  \( a_x \) and \( b_x \) be linear forms such that \( \{a_x,b_x\} \)  is linearly independent. The following are equivalent:
 \begin{enumerate}
 \item \( \J3[f_{xxx},a_x^2b_x,u_x^3]=0 \)
 
 \item \( f_{xxx}=  a_x^2 d^{\phantom{2}}_x + b^{\phantom{2}}_0 b_x^3\) for some linear form \( d_x \) and \( b_0\in \C \).
 \end{enumerate}
 \end{lemma}
 
 \begin{proof}
If  \( f_{xxx} \) satisfies (2),  then (1) follows from \eqref{eqJ} and \eqref{eqL}.

Let \( c_x \) be a linear form such that  \( \{a_x,b_x,c_x\} \)  is linearly independent. Given \( \J3[f_{xxx},a_x^2b_x,u_x^3]=0 \), the identity
  \begin{align*}
 \con_{ux}^3[&\J3[f_{xxx}, a_x^2 b_x, u_x^3]\, a_x b_x^2 c_x^3] =\\
&  12 \Big(-\J3[f_{xxx}, a_x^3, b_x^3]\, c_x^3 + 3 \J3[f_{xxx}, b_x^3, a_x^2 c_x]\, c_x^2 a_x  \\
&\qquad- 3 \J3[f_{xxx}, a_x^3, c_x^2 b_x]\, b_x^2 c_x
- 3 \J3[f_{xxx}, c_x^3, a_x^2 b_x]\, b_x^2 a_x  \\
&\qquad- 6 \J3[f_{xxx}, a_x^3, b_x^2 c_x]\, c_x^2 b_x - 18 \J3[f_{xxx}, a_x b_x^2, a_x c_x^2]\, a_x b_x c_x\Big) 
 \end{align*}
and the linear independence of \( \B_3 \) imply that six of the terms on the right side of \eqref{equ33} are zero. The remaining terms are either multiples of \(a_x^2\) or multiples of \(b_x^3\), so \( f_{xxx} \) can be written in the claimed form.
 \end{proof}
 
\section{Reducibility of Quadratic Forms}

This section serves to put well-known results about the reducibility of quadratic forms into the language of transvectants--- the language we will use for cubic forms. See \cite{kronenthal} for a general discussion of quadratic forms, including forms in more than three variables. 

For a quadratic form \( f_{xx} \) (with coefficients labelled as in \eqref{eq00}), there are only two important concomitants:
\[ \J2[f_{xx},f_{xx},f_{xx}]=12 (4 f_{11} f_{22} f_{33}+ f_{12}  f_{23}f_{13}   -f_{23}^2 f_{11}  -f_{13}^2 f_{22}- f_{12}^2 f_{33} ) \] 
and
\begin{align*}
\J2[f_{xx},f_{xx},u_x^2]&=4 ( 4 f_{22} f_{33}-f_{23}^2 ) u_1^2 +4 ( 4 f_{11} f_{33}-f_{13}^2 ) u_2^2\\
&\qquad+ 4 ( 4 f_{11} f_{22}-f_{12}^2) u_3^2+ 8 (f_{13} f_{23} - 2 f_{12} f_{33}) u_1 u_2 \\
&\qquad  +  8 (f_{12} f_{13} - 2 f_{11} f_{23}) u_2 u_3+ 8 ( f_{12} f_{23}-2 f_{13} f_{22}) u_1 u_3 
\end{align*}
The first of these concomitants is the Hessian \eqref{eq388} of \( f_{xx} \), but it is also frequently called the discriminant of \( f_{xx} \). The \textbf{discriminant} of an arbitrary ternary form \( f \),  is, by definition, a concomitant of \( f \) that is zero if and only if \( f \) is \textbf{singular,} that is, all partial derivatives \( (\partial_1 f, \partial_2 f,\partial_3 f)\) are zero for some \( \xx\in \C \), not all zero.  For quadratic forms, it turns out that the discriminant, the Hessian and \( \J2[f_{xx},f_{xx},f_{xx}] \) coincide (up to multiplication by constants)~\cite{kronenthal}.  

\new We say, somewhat imprecisely,  that \( f_{xx} \) is a sum of three squares if 
 \[ f_{xx}=a_0 a_x^2+ b_0 b_x^2 +c_0 c_x^2 \]
for some linear forms \( a_x \), \( b_x \) and \( c_x \), and \( a_0,b_0, c_0\in \C \). Similarly, we say  \( f_{xx} \) is sum of two squares if \( f_{xx}=a_0 a_x^2+ b_0 b_x^2 \) and \( f_{xx} \) is a square if \( f_{xx}=a_0 a_x^2 \).

Supposing that \( f_{xx} \) is a sum of three squares as above, from \eqref{eqJ}, \eqref{eqK} and the multilinearity of \( \JJ2\), we get 
\begin{equation}\label{eqw51}
\begin{aligned}
\J2[f_{xx},f_{xx},u_x^2]&=16\left( a_0 b_0 [abu]^2+ b_0 c_0 [bcu]^2+  c_0 a_0 [cau]^2\right)\\[5pt]
\J2[f_{xx}, f_{xx}, f_{xx}] &= 48 a_0 b_0 c_0 [abc]^2
\end{aligned}
\end{equation}
If \( f_{xx} \) is a sum of two squares, then we can assume that \( c_0=0 \) in the above, and so  \( \J2[f_{xx}, f_{xx}, f_{xx}]=0 \).  If \( f_{xx} \) is a square, then we can assume that \( b_0=c_0=0 \) in the above, and so  \( \J2[f_{xx}, f_{xx}, u_x^2]=0 \).  The converses of these statements are true as is proven in \lref{lem10} below.

%

For the  proof, we note that, if \( f_{xx} \) is any quadratic form  and \( \yy\in \C \) are fixed, then  \( g_{xx}=4f_{yy}f_{xx}-f_{xy}^2 \) is a quadratic form in \( \xx \)  that satisfies
\begin{gather}
4g_{xx}=\left\Vert \J2[f_{xx},f_{xx},u_x^2]\right\Vert_{u\mapsto [xy]} \label{eqq11}\\
3 \J2[g_{xx}, g_{xx}, u_x^2] = 16\J2[f_{xx}, f_{xx}, f_{xx}] f_{yy} u_y^2 \label{eqq12}\\
\J2[g_{xx},g_{xx},g_{xx}] = 0\label{eqq13}
\end{gather}


\begin{lemma}\label{lem10}
Let \( f_{xx} \) be a quadratic  form. 
\begin{enumerate}
\item \( \J2[f_{xx},f_{xx},u_x^2]=0 \)  if and only if  \( f_{xx}=a_0a_x^2 \) for some linear form \( a_x \) and \( a_0\in \C \).

\item  \( \J2[f_{xx},f_{xx},f_{xx}]=0 \) if and only if  \( f_{xx}=a_0a_x^2+ b_0 b_x^2 \) for some linear forms \( a_x \), \( b_x \) and \( a_0,b_0 \in \C \).

\item \( f_{xx}=a_0 a_x^2+ b_0 b_x^2 +c_0 c_x^2\)  for some linear forms \( a_x \), \( b_x \), \( c_x \), and \( a_0,b_0, c_0\in \C \).
\end{enumerate}
\end{lemma}

\begin{proof}

In (1) and (2), if \( f_{xx} \) has the form on the right, then the transvectant on the left is zero because of  \eqref{eqJ}, \eqref{eqK} and the multilinearity of transvectants.

Fix \( \yy\in \C \) so that \( f_{yy}\in \C \) is nonzero and set  \( g_{xx}=4f_{yy}f_{xx}-f_{xy}^2 \). 

(1) Suppose that  \( \J2[f_{xx},f_{xx},u_x^2]=0 \). Because of the identity \eqref{eqq11}, we have  \( g_{xx}=0 \), that is, \( 4f_{yy}f_{xx}-f_{xy}^2=0 \).  Solving this equation for \( f_{xx} \) we get \( f_{xx}=a_0a_x^2 \)  with \( a_x=f_{xy} \). 

(2) Suppose that \( \J2[f_{xx}, f_{xx}, f_{xx}]=0 \). Because of (1) of this lemma and the identity \eqref{eqq12} we have \( g_{xx}= a_0a_x^2 \) for some form \( a_x \) and \( a_0\in \C \), that is, \( 4 f_{xx} f_{yy} - f_{xy}^2 = a_0a_x^2 \). Solving this equation for \( f_{xx} \) we get  \( f_{xx}=a_0a_x^2+ b_0 b_x^2 \) with \( b_x=f_{xy} \) and the same \( a_x \), but different \( a_0 \). 
 
 (3) Because of (2) of this lemma and the identity \eqref{eqq13} we have \( g_{xx}= a_0a_x^2 +b_0 b_x^2\) for some linear forms \( a_x \), \( b_x \) and \( a_0,b_0 \in \C \), that is, \( 4 f_{xx} f_{yy} - f_{xy}^2 =a_0a_x^2 +b_0 b_x^2 \). Solving this equation for \( f_{xx} \) we get  \( f_{xx}=a_0a_x^2+ b_0 b_x^2+c_0c_x^2  \) with \( c_x=f_{xy} \) and the same \( a_x \) and \( b_x \), but different \( a_0 \) and \( b_0 \). 
\end{proof}


The expression \( \J2[f_{xx},f_{xx},u_x^2] \) has six coefficients, and so, by \lref{lem10}(1), to determine if  \( f_{xx} \) is a square, six corresponding equations have to be checked. Alternatively, from the proof of \upitem1, we see that \( f_{xx} \) is a square if and only if  \( g_{xx}=4f_{yy}f_{xx}-f_{xy}^2 \) is zero for some \( \yy\in \C \) such that \( f_{yy} \) is nonzero. With  \( \yy \) fixed,  \( g_{xx} \) is a quadratic form in \( \xx \) so also has six coefficients to check. 

In fact, generically, it suffices to check only three equations.  For example, if \( f_{33} \) is nonzero  we can set \( (\yy)=(0,0,1) \) so that \( f_{yy}=f_{33}\neq 0 \). Then 
\[  g_{xx}=( 4 f_{11} f_{33}-f_{13}^2) x_1^2 + 
 2 ( 2 f_{12} f_{33}-f_{13} f_{23} ) x_1 x_2 + ( 4 f_{22} f_{33}-f_{23}^2 ) x_2^2 \] 
 So \( f_{xx} \) is a square if and only if \(  4 f_{11} f_{33}-f_{13}^2= 2 f_{12} f_{33}-f_{13} f_{23} = 4f_{22} f_{33}-f_{23}^2 =0 \) 
 and in this circumstance \[  f_{xx}=\frac1{4f_{33}}(f_{13} x_1 + f_{23} x_2 + 2 f_{33} x_3)^2 \] 

 \begin{lemma}\label{lemw78}
 A quadratic form \( f_{xx} \) is reducible if and only if \( \J2[f_{xx},f_{xx},f_{xx}]=0 \).
 \end{lemma}
 \begin{proof}
If \( f_{xx} \) is a sum of two squares, that is,   \( f_{xx}=a_0a_x^2+ b_0 b_x^2 \) for some linear forms \( a_x \), \( b_x \) and \( a_0,b_0 \in \C \), then \( f_{xx}=(\alpha_0 a_x+i \beta_0 b_x)(\alpha_0 a_x-i \beta_0 b_x) \)
where \( \alpha_0^2=a_0 \) and \( \beta_0^2=b_0 \), and so \( f_{xx} \) is reducible. Conversely, if \( f_{xx} \) is reducible, that is, \( f_{xx}=a_xb_x \) for some linear forms \( a_x \), \( b_x \),  then \( 4f_{xx}=(a_x+b_x)^2-(a_x-b_x)^2 \) so \( f_{xx} \) is a sum of two squares.

The claim now follows from \lref{lem10}(2).
 \end{proof}

\begin{lemma}
A quadratic form \( f_{xx} \) is irreducible if and only if there are linearly independent linear forms  \( \{a_x,b_x,c_x\} \) and nonzero  \( a_0,b_0, c_0\in \C \)
such that  \begin{equation}\label{eqw43}
 f_{xx}=a_0 a_x^2+ b_0 b_x^2 +c_0 c_x^2.
\end{equation}
\end{lemma}

\begin{proof}
Because of \lref{lem10}, any quadratic form can be written as in \eqref{eqw43}  and  in that case, by \eqref{eqw51},   \( \J2[f_{xx},f_{xx},f_{xx}]=48 a_0b_0c_0\,[abc]^2 \). By \lref{lemw78}, \( f_{xx} \) is irreducible if and only if \( a_0b_0c_0\,[abc]\neq 0 \), if  and only if \( a_0,b_0, c_0\in \C \) are nonzero and  \( \{a_x,b_x,c_x\} \) is linearly independent. 
\end{proof}

%

The following lemma is not directly related to the reducibility question, but is nonetheless worth mentioning. It  makes clear the difference between the hypotheses  \( \J2[f_{xx},f_{xx},u_x^2]=0  \) (\lref{lem10}(1)) and  \( \J2[f_{xx},f_{xx},a_x^2]=0  \). In the first case, the coefficients of \( u_x \) are variables, so  \( \J2[f_{xx},f_{xx},u_x^2] \) is a quadratic form in \( \uu \) with six coefficients. The assumption that \( \J2[f_{xx},f_{xx},u_x^2]=0  \) means that all six coefficients are zero. In contrast, the coefficients of \( a_x \) are fixed, so \( \J2[f_{xx},f_{xx},a_x^2]  \) is a single number and  the hypothesis  \( \J2[f_{xx},f_{xx},a_x^2] =0 \), in the following lemma,  is much weaker assumption.

 \begin{lemma}\label{lem93}
 Let \( f_{xx} \) be a quadratic form and let \( a_x \) be a linear form.
 Then \( \J2[f_{xx},f_{xx},a_x^2]=0 \) if and only if  \( f_{xx}=a_x b_x+ c_0 c_x^2 \) for some linear forms \( b_x \) and \( c_x \) and  \( c_0\in \C \). 
 \end{lemma}

\begin{proof}
 If \( f_{xx}=a_x b_x+ c_0 c_x^2 \), then  \( \J2[f_{xx},f_{xx},a_x^2]=0 \) follows by  direct calculation, or from \eqref{eqJ}.

 For the converse, fix \( \yy\in \C \) so that \( a_y\in \C \) is nonzero and set \[ g_{xx}= a_x^2 f_{yy} -   a_x a_y f_{xy}+  a_y^2 f_{xx} .\] Then \( g_{xx} \) is a quadratic form in \( \xx \) that satisfies the identity
 \[\J2[g_{xx},g_{xx},u_x^2]= a_y^2 u_y^2 \J2[f_{xx}, f_{xx}, a_x^2] \] 
 Since \( \J2[f_{xx}, f_{xx}, a_x^2]=0 \), we have \(  \J2[g_{xx},g_{xx},u_x^2] =0 \) and so, because of  \lref{lem10}(1),  \( g_{xx}= c_0c_x^2 \) for some  linear form \( c_x \) and \( c_0\in \C \). The resulting equation, \( c_0c_x^2=a_x^2 f_{yy} -   a_x a_y f_{xy}+  a_y^2 f_{xx}  \) can be solved for \( f_{xx} \) putting it into the claimed form. 
\end{proof}

If  \( f_{xx}=a_x b_x+ c_0 c_x^2 \), then \( \J2[f_{xx},f_{xx},f_{xx}] =-12c_0 [abc]^2 \), so by \lref{lem10}(2), \( f_{xx} \) is irreducible if and only if \(  c_0\neq 0 \) and \( \{a_x,b_x,c_x\} \) is linearly independent. When this happens, the line \( a_x=0 \) is tangent to the curve \( f_{xx}=0 \) at the intersection point of the lines \( a_x=0 \) and \( c_x=0 \).

%

\section{The Concomitants of a Cubic Form}\label{cubic}

We now introduce the most important concomitants that can be constructed from a single cubic form \( f_{xxx} \).
\begin{equation}\label{eq837}
\begin{aligned}
\theta_{uuxx}&=\frac14\J2[f_{xxx},f_{xxx},u_x^2]&
\Delta_{xxx}&=\frac 1{12}\J2[f_{xxx},f_{xxx},f_{xxx}]\\
S_{uuu}&=\frac1{576}\J4[f_{xxx} u_x,f_{xxx}u_x,f_{xxx}u_x]&
S&=\Vert S_{uuu} \Vert_{u^3\mapsto f}\\
T_{uuu}&=\frac1{576}\J4[f_{xxx} u_x,f_{xxx}u_x,\Delta_{xxx}u_x]&
T&=\Vert T_{uuu} \Vert_{u^3\mapsto f}\\
\Pi_{u4x} &=\frac1{12} \J[\Delta_{xxx}, f_{xxx}, u_x] &
\Gamma_{4ux} &=\frac1{432}\J3[ \Pi_{u4x},f_{xxx},u_x^3]
\end{aligned}
\end{equation}
\[ F_{6u} =\frac1{3072} \J2[\theta_{uuxx},\theta_{uuxx},u_x^2] \]
The rational coefficients in these definitions serve merely to eliminate integer factors that would otherwise be common to all terms.  For  \( \Pi_{u4x}\), \( \Gamma_{4ux}  \) and \(  F_{6u} \) we have extended the notational convention in an obvious way because \( \Pi_{uxxxx}\), \( \Gamma_{uuuux} \)  and  \(  F_{uuuuuu} \) are too cumbersome. The symbols chosen for the concomitants largely follow the 19th century German literature, except for \( \Pi_{u4x} \) and \( \Gamma_{4ux} \).

Of course, \( \Delta_{xxx} \) is the Hessian \eqref{eq388} of \( f_{xxx} \) (we ignore the rational coefficient), first investigated by Hesse in 1844 \cite{hesse1844a,hesse1844b,hesse1852}.  The concomitants \( \theta_{uuxx} \) and \(  F_{6u}  \) were first calculated by Cayley~\cite{cayley1846} in 1846. Aronhold \cite{aronhold1849} introduced the invariants \( S \) and \( T \) in 1849. They were written out explicitly for the first time by Salmon \cite{salmon1851} in 1851. 
Cayley collected formulas for all these concomitants, except \( \Pi_{u4x} \) and \(  \Gamma_{4ux}\), in his \textit{Third Memoir on Quantics} \cite{cayley1856} in 1856. The concomitant \( S_{uuu} \) is called the \textbf{Cayleyan} by many authors (excluding Cayley, who called  \( S_{uuu} \)  the Pippian and \( T_{uuu} \) the Quippian).  


For each of the concomitants, the table below gives its degree with respect to the variables \( \xx \), the variables \( \uu \), and the coefficients of \( f_{xxx} \). The bottom row gives the number of terms in the concomitant when fully expanded.
\[ \begin{array}{c|c|c|c|c|c|c|c|c|c|c|c|}
 & f_{xxx} & \Delta_{xxx} & \theta_{uuxx} & S_{uuu} & T_{uuu} & S & T  &\Pi_{u4x}  & \Gamma_{4ux}&  F_{6u} \\
 \hline
x & 3 & 3 & 2 & 0 & 0 & 0 & 0  & 4&1 & 0 \\
u & 0 & 0 & 2 & 3& 3& 0 & 0  & 1& 4& 6\\
f & 1 & 3 & 2 & 3 & 5 & 4 & 6  & 4&5& 4\\\hline
\text{terms} & 10 & 73 & 84 & 82 & 448 & 25 & 103  & 576 &1314 & 418
\end{array}
 \] 
As well as the determinantal expressions for \( \Delta_{xxx} \) and \( \theta_{uuxx} \) in \eqref{eq388} and \eqref{eq382}, there is a useful expression for \( S_{uuu} \) due to Aronhold \cite[p.~189]{aronhold1858}:
\[ S_{uuu}=\begin{vmatrix}
3 f_{111}&  f_{112}&  f_{113}&  u_1&  0&  0 \\ 
2 f_{112}&  2 f_{122}&  f_{123}&  u_2&  u_1&  0 \\
 f_{122}&  3 f_{222}&  f_{223}&  0&  u_2&  0 \\
 2 f_{113}&  f_{123}&  2 f_{133}&  u_3&  0& u_1 \\
 f_{123}&  2 f_{223}&  2 f_{233}&  0&  u_3&  u_2 \\
 f_{133}&  f_{233}&  3 f_{333}&  0&  0&  u_3
\end{vmatrix}
 \] 
%
From the expression for \( S_{uuu} \) it is easy to calculate \( S \) and \( T \) using the identities 
\[ S=\Vert S_{uuu}\Vert_{u^3\mapsto f}\quad T=\Vert S_{uuu}\Vert_{u^3\mapsto \Delta}  \]
where the substitution \( u^3\mapsto \Delta \) replaces the coefficients of \( u_x^3 \) by the corresponding coefficients of \( \Delta_{xxx} \). For example,
 \begin{align*}
S&=f_{123}^4 - 8 f_{122} f_{123}^2 f_{133} + 16 f_{122}^2 f_{133}^2 +  24 f_{113} f_{123} f_{133} f_{222} - 48 f_{112} f_{133}^2 f_{222}\\
&\quad  - 8 f_{113} f_{123}^2 f_{223} - 16 f_{113} f_{122} f_{133} f_{223} + 24 f_{112} f_{123} f_{133} f_{223} + 16 f_{113}^2 f_{223}^2\\
&\quad  -  48 f_{111} f_{133} f_{223}^2 + 24 f_{113} f_{122} f_{123} f_{233} - 8 f_{112} f_{123}^2 f_{233} - 16 f_{112} f_{122} f_{133} f_{233}\\
&\quad  - 48 f_{113}^2 f_{222} f_{233} +  144 f_{111} f_{133} f_{222} f_{233} - 16 f_{112} f_{113} f_{223} f_{233}\\
&\quad  +  24 f_{111} f_{123} f_{223} f_{233} + 16 f_{112}^2 f_{233}^2 - 48 f_{111} f_{122} f_{233}^2 \\
&\quad -  48 f_{113} f_{122}^2 f_{333} + 24 f_{112} f_{122} f_{123} f_{333} +  144 f_{112} f_{113} f_{222} f_{333} \\
&\quad - 216 f_{111} f_{123} f_{222} f_{333} -  48 f_{112}^2 f_{223} f_{333} + 144 f_{111} f_{122} f_{223} f_{333}
\end{align*}
Formulas for \( S \) and \( T \) can be found in \cite[Prop.~4.4.7, Ex.~4.5.3]{sturmfels}.

There are many other expressions that could be used to calculate, or define, these concomitants, for example,
\begin{align*}
\Delta_{xxx}&=\frac1{2}\J[ \partial_1 f_{xxx}, \partial_2 f_{xxx}, \partial_3 f_{xxx}] &
S&= \frac1{576}\con_{ux}^4[\theta_{uuxx}^2]\\
S_{uuu}&=-\frac1{24}\J3[u_x \partial_1 f_{xxx},u_x \partial_2 f_{xxx},u_x \partial_3 f_{xxx}] &
T&= -\frac1{276\,480} \con_{ux}^6[\theta_{uuxx}^3]\\
S_{uuu}&=\frac38 \con_{ux}[\J2[\theta_{uuxx},f_{xxx},u_x^2]] &
\llap{$T_{uuu} =\dfrac38$} &\con_{ux}[\J2[\theta_{uuxx},\Delta_{xxx},u_x^2]] \\
\Pi_{u4x}&=\frac1{576 } \J3[f_{xxx} , f_{xxx}^2 , f_{xxx}  u_x] &
\llap{$ \Gamma_{4ux}=  \dfrac{1}{24}$} &\left(2 T_{uuu} u_x - \con_{ux}[S_{uuu} \theta_{uuxx}]\right)\\
F_{6u}&= \frac1{1728} \J3[\J[ \theta_{uuxx}, f_{xxx}, u_x], f_{xxx}, u_x^3]
\end{align*}

\new  Of course, the reason for introducing all these concomitants of \( f_{xxx} \) is that they determine whether \( f_{xxx} \) is completely reducible and how   \( f_{xxx} \) factors when it is completely reducible. Though the proofs will come later, it is useful  to collect these relationships in one theorem.
\begin{theorem}
Let \( f_{xxx} \) be a ternary cubic form. 
\begin{enumerate}
\item \( \theta_{uuxx}=0 \) if and only if  \( f_{xxx}=a_0a_x^3 \) for some linear form \( a_x \) and \( a_0 \in \C \).

\item \(  F_{6u} =0 \) if and only if \( f_{xxx}=a_x b_x^2 \) for some linear forms \( a_x \) and \( b_x \).

\item  \( \Delta_{xxx}=0 \) if and only if \( S_{uuu}=0 \), if and only if \( f_{xxx}=a_xb_xc_x \) for some linearly dependent linear forms \( a_x \), \( b_x \) and \( c_x \).

\item  \( \Pi_{u4x} =0 \) if and only if   \( \Gamma_{4ux}=0 \), if and only if \( f_{xxx}=a_xb_xc_x \) for some linear forms \( a_x \), \( b_x \) and \( c_x \).
\end{enumerate}
\end{theorem}

\begin{proof}
Lemmas~\ref{lem50x}, \ref{lem381}, \ref{lemu93} and \tref{thmq111}, respectively.
\end{proof}

%


\new We now consider the values of these concomitants when \( f_{xxx} \) is completely reducible. If \( f_{xxx}=a_xb_xc_x \) for some linear forms \( a_x \), \( b_x \) and \( c_x \), then 
\begin{align}
\Delta_{xxx}&=[abc]^2\, a_xb_xc_x\label{eq370}\\[3pt]
\theta_{uuxx}&= \rlap{\( [abc]^2 u_x^2-2 \left( [bcu]^2 a_x^2 + [cau]^2 b_x^2+[abu]^2 c_x^2 \right) \)} \label{eq390}\\[3pt]
S_{uuu}&=[abc]\,[abu]\,[bcu]\,[cau] & S&=[abc] ^4 \notag \\[3pt]
T_{uuu}&=[abc]^3\,[abu]\,[bcu]\,[cau] & T&=[abc] ^6 \notag \\[3pt]
 F_{6u} &=[abu]^2[bcu]^2[cau]^2 &  \Pi_{u4x} &= \Gamma_{4ux}=0\label{eq376}
\end{align}

As seen in \eqref{eq370} and pointed out in the introduction, if \( f_{xxx} \) is completely reducible, then \( \Delta_{xxx} \) is a multiple of \( f_{xxx} \), specifically, \( \Delta_{xxx}= [abc]^2 f_{xxx}\).  It is easy to find other necessary conditions for the complete reducibility of \( f_{xxx} \). For example,
\begin{equation}\label{eqy5}
\begin{gathered}
S^3-T^2 =0 \qquad S \Delta_{xxx}-T f_{xxx} =0 \qquad  S T_{uuu}-T S_{uuu}=0 \\
  T \Delta_{xxx}-S^2 f_{xxx} =0 \qquad  S S_{uuu}^2 -  T  F_{6u} =0\\
  \J2[\Delta_{xxx}, \Delta_{xxx}, u_x^2] - 4 S \theta_{uuxx}=0 \qquad 
 S   F_{6u} -S_{uuu}T_{uuu}=0 
\end{gathered}
\end{equation}


The hope that any of these conditions are sufficient for the complete reducibility of \( f_{xxx} \) is spoiled by the following example.


 \begin{example}
 If \( f_{xxx}=x_1(x_1 x_2+x_3^2) \),  then 
 \[  \Delta_{xxx}=-4 x_1^3\qquad S=T=0 \qquad T_{uuu}=0\qquad \J2[\Delta_{xxx}, \Delta_{xxx}, u_x^2]=0.\] Since \( \Delta_{xxx} \) is not a multiple of \( f_{xxx} \), \( f_{xxx} \) is not completely reducible even though all the equations in \eqref{eqy5} hold for this particular  form.
 \end{example}
 
 It turns out that the similar looking conditions  
\( F_{6u} \Delta_{xxx} - S_{uuu}^2 f_{xxx}=0\) and \(\Delta_{xxx}^2-S f_{xxx}^2=0  \)  are each necessary and sufficient for the complete reducibility of \( f_{xxx} \) (\lref{lemq120}), 
and, in fact, several of the conditions in \eqref{eqy5} are necessary and sufficient for the complete reducibility of \( f_{xxx} \) so long as \( S \) is nonzero (\lref{lemq121}). 

\new We now consider the concomitant \(\theta_{uuxx} \).

\begin{lemma}\label{lem50x}
For a cubic form \( f_{xxx} \) the following are equivalent:
\begin{enumerate}
 \item \(\theta_{uuxx}=0 \)
   
 \item \( 27 f_{yyy}^2 f_{xxx} - f_{xyy}^3=0  \)
 
 \item \( f_{xxx}=a_0a_x^3 \) for some linear form \( a_x \) and \( a_0 \in \C \).
\end{enumerate}
\end{lemma}

\begin{proof}
 If (1) holds, then (2) follows from the identity
\[ 4 (27 f_{yyy}^2 f_{xxx} - f_{xyy}^3) = \Vert 3 f_{yyy} \theta_{uuxy} + f_{xyy} \theta_{uuyy} \Vert_{u\mapsto[xy]}  \] 
If (2) holds, then we can fix \( \yy \in \C \) so that \( f_{yyy}\in \C \) is nonzero. Then \[ f_{xxx}=\frac1{27f_{yyy}}  f_{xyy}^3\]
showing that \( f_{xxx}=a_0a_x^3\)  for some linear form \( a_x \) and \( a_0 \in \C \) as claimed in (3).

If (3) holds, then (1) follows from the definition of \( \theta_{uuxx}  \) in \eqref{eq837} and  \eqref{eqJ}. 
\end{proof}

It seems that the above result was first published in 1871 by Gundelfinger~\cite[p.~571] {gundelfinger1871b}, though it is hard to imagine that it had not been noticed earlier.

It is worth noting  that \( \theta_{uuxx} \) is a quadratic form in \( \xx \), and also in \( \uu \) so has potentially \( 6\cdot 6= 36 \) coefficients that have to be checked to confirm that \( \theta_{uuxx}=0 \). In fact, these coefficients are not linearly independent because the identity \( \con_{ux}[ \theta_{uuxx}]=0 \)  gives \( 9 \)  linear equations satisfied by the coefficients.  So it suffices to check only 27 coefficients to confirm that \( \theta_{uuxx}=0 \).


The expression \( 27 f_{yyy}^2 f_{xxx} - f_{xyy}^3  \) is degree \( 6 \) in \( y_1 \), \( y_2 \) and \( y_3 \), and  degree \( 3 \) in \( x_1 \), \( x_2 \) and \( x_3 \) so has potentially \( 28\cdot10=280 \) coefficients that have to be checked to determine if  \( 27 f_{yyy}^2 f_{xxx} - f_{xyy}^3=0  \). (In fact, there are \( 262 \) nonzero coefficients.)

But the condition \( 27 f_{yyy}^2 f_{xxx} - f_{xyy}^3=0  \) has the great advantage that we do not have to confirm it for all \( \yy \in \C \).  The  proof of \lref{lem50x}  shows  that, if we fix \( \yy \in \C \) such that \( f_{yyy}\neq 0 \), then \( f_{xxx} \) is a cube if and only if  the ten coefficients of \( 27 f_{yyy}^2 f_{xxx} - f_{xyy}^3  \) are zero.

In fact,  generically, seven coefficients suffice. For example, if \( f_{333}\neq 0 \), then we can choose \( (\yy)=(0,0,1) \) so that  \( f_{yyy}=f_{333} \) is nonzero. Then 
\begin{align*}
27 f_{yyy}^2 f_{xxx} - f_{xyy}^3 &= ( 27 f_{111} f_{333}^2-f_{133}^3) x_1^3 + ( 27 f_{222} f_{333}^2-f_{233}^3) x_2^3 \\
&\quad+ 3 ( 9 f_{122} f_{333}^2-f_{133} f_{233}^2 ) x_1 x_2^2 + 3 ( 9 f_{112} f_{333}^2-f_{133}^2 f_{233} ) x_1^2 x_2   \\
&\quad+ 9 f_{333} ( 3 f_{113} f_{333}-f_{133}^2 ) x_1^2 x_3 +  9 f_{333} ( 3 f_{223} f_{333}-f_{233}^2 ) x_2^2 x_3 \\
&\quad + 9 f_{333} ( 3 f_{123} f_{333}-2 f_{133} f_{233} ) x_1 x_2 x_3,
\end{align*}
So \( f_{xxx} \) is a cube if and only if the seven coefficients in the above expression are zero. When this factorization occurs, we get \( 27 f_{333}^2f_{xxx}=(f_{133} x_1 + f_{233} x_2 + 3 f_{333} x_3)^3 \).

\section{The Hessian of a Cubic Form}\label{hessian}

The most important result to be proved in this section is that,  for a cubic form \( f_{xxx} \),  its Hessian \( \Delta_{xxx} \) is  zero if and only if  \( f_{xxx}=a_xb_xc_x \) for some linearly dependent forms  \( \{a_x,b_x,c_x\}\). It seems that  there are two ways to prove this, either by considering  \( \theta_{uuxx} \)  as a quadratic form in \( \xx \) for some fixed \( \uu\in \C \) (\lref{lemu93}), or as a quadratic form in \( \uu \) for some fixed \( \xx\in \C \) (\lref{lemg91}).  In the following lemma we choose the first method.
\begin{lemma}\label{lemu93}
For a cubic form  \( f_{xxx} \), the following are equivalent:
\begin{enumerate}
\item \( \Delta_{xxx}=0 \)

\item \( S_{uuu}=0 \)

\item \( f_{xxx}=a_0 a_x^3+b_0 b_x^3 \) for some \( a_0,b_0\in \C \) and linear forms \( a_x \) and \( b_x \), or \( f_{xxx}=a_x^2 b_x \) for some  linear forms \( a_x \) and \( b_x \). 

\item \( f_{xxx}=a_xb_xc_x \)  for some linearly dependent forms  \( \{a_x,b_x,c_x\}\).
\end{enumerate}
\end{lemma}

\begin{proof}

\noindent(1)\( \Leftrightarrow \)(2): The equivalence of the conditions \( \Delta_{xxx}=0 \) and  \( S_{uuu}=0 \) follows from the identities
\begin{align*}
768 \,S_{uuu}^2&=\J2[\J2[\Delta_{xxx}, f_{xxx}, u_x^2], \theta_{uuxx}, u_x^2] \\
18\, \Delta_{xxx}^2 &= 111 \con_{ux}^3[S_{uuu} f_{xxx}] f_{xxx}^2- 4 \con_{ux}^3[S_{uuu} f_{xxx}^3] 
\end{align*}

\noindent(1,2)\( \Rightarrow \)(3): If \( \theta_{uuxx}=0 \), then, by \lref{lem50x}, \( f_{xxx}=a_0 a_x^3 \) for some \( a_0\in \C \) and linear form \( a_x \), and so \upitem3 holds (in two different ways).  Otherwise, fix \( \vv\in \C \) so that \( \theta_{vvxx} \) is a nonzero quadratic form in \( \xx \). The identity
\[ 27 \J2[\theta_{vvxx}, \theta_{vvxx}, \theta_{vvxx}] = 16 S_{vvv}^2 \] and \lref{lemw78} imply that  \( \theta_{vvxx}=a_xb_x \) for some nonzero linear forms  \( a_x \) and \( b_x \).
The identity
\[ 27 \J2[f_{xxx}, \theta_{vvxx}, u_x^2] =8 S_{uuv} v_x - 4 S_{uvv} u_x + 144 \left\Vert \Delta_{xyz} \right\Vert_{\substack{y\mapsto[u,v]\\z\mapsto[u,v]}}
 \]then implies that \(  \J2[f_{xxx}, a_xb_x, u_x^2]=0\).  There are now two cases: If \( \{a_x,b_x\} \) is linearly dependent, then \(  \J2[f_{xxx}, a_x^2, u_x^2]=0\) and so by \lref{lemu91}(C2),  \(  a_x^2 \) divides \( f_{xxx} \). In the remaining case, \( \{a_x,b_x\} \) is linearly independent and so by \lref{lemu46}(C), \( f_{xxx}=a_0 a_x^3+b_0 b_x^3 \) for some \( a_0,b_0\in \C \). 
 
\noindent(3)\( \Rightarrow \)(4): Either \( f_{xxx} =a_x^2 b_x\) or 
\[ f_{xxx}=a_0 a_x^3+b_0 b_x^3=(\alpha_0 a_x+\beta_0 b_x)(\omega\alpha_0  a_x+ \omega^2 \beta_0 b_x) (\omega^2\alpha_0  a_x+\omega \beta_0  b_x)  \] where \( \alpha_0^3=a_0 \), \( \beta_0^3=b_0 \) and \( \omega=e^{2\pi i/3} \).
Either way, \( f_{xxx} \) is completely reducible with  factors that are linearly dependent so \upitem4  holds.

\noindent(4)\( \Rightarrow \)(1): This follows immediately from \eqref{eq370}.
\end{proof}

Sylvester \cite[p.~187]{sylvester1852} seems to be the first (and possibly only) person to claim that the conditions \upitem1 and \upitem2  from this lemma are equivalent. 

\new The concomitant \(  F_{6u} \) distinguishes the two cases in \lref{lemu93}\upitem3.

\begin{lemma}\label{lem381}
For a cubic form  \( f_{xxx} \), the following are equivalent:
\begin{enumerate}
\item   \(  F_{6u} =0 \)

\item  \( f_{xxx}=a_x b_x^2 \) for some linear forms \( a_x \)  and \( b_x \).
\end{enumerate}
\end{lemma}

\begin{proof}
Suppose that \(  F_{6u} =0 \). Because of the identity
\(  S_{uuu}^2 =36 \con_{ux}^2[ F_{6u}  \theta_{uuxx}] \)
we have \( S_{uuu}=0 \), and then, by \lref{lemu93}, \( f_{xxx}=a_xb_x c_x \) for some nonzero linear forms \( a_x \), \( b_x \) and \( c_x \). 
Because of \eqref{eq376},  \(  F_{6u} =0 \) now implies that \( [abu]=0 \), \([bcu]=0\) or \( [cau]=0 \). Thus at least one pair from \( \{a_x,b_x,c_x\} \) is linearly dependent. After some relabeling and scaling, \( f_{xxx} \) can be written in the claimed form. 

If  \( f_{xxx}=a_x b_x^2 \) for some linear forms \( a_x \)  and \( b_x \), then \(  F_{6u}=0 \) follows immediately from  \eqref{eq376}.
\end{proof}

 Perhaps worth noting is that, if \(   f_{xxx}=a_0 a_x^3+b_0 b_x^3 \), then 
 \begin{equation}\label{eqq133}
   F_{6u} =-27 a_0^2 b_0^2\, [abu]^6
 \end{equation}
If \(  F_{6u} =0 \), then either \( a_0=0 \), \(b_0=0\) or \( [abu]=0 \). In all these cases, \( f_{xxx} \) is a cube. In other words, if \( f_{xxx} \)  can be written both as \( f_{xxx}=a_0 a_x^3+b_0 b_x^3 \) and as \( f_{xxx}=a_x b_x^2 \) for some possibly different linear forms \( a_x \)  and \( b_x \), then \( f_{xxx} \) is a cube.

\begin{lemma}\label{lemq27}
For a cubic form  \( f_{xxx} \), the following are equivalent:
\begin{enumerate}
\item \( \Delta_{xxx}=0 \) and  \(  F_{6u} \neq 0\).

\item \( f_{xxx}=a_0 a_x^3+b_0 b_x^3 \) for some nonzero \( a_0,b_0\in \C \) and linear forms \( a_x \) and \( b_x \) such that \( \{a_x,b_x\} \) is linearly independent.

\item \( f_{xxx}=a_xb_xc_x \)  for some linearly dependent, but pairwise independent, forms  \( \{a_x,b_x,c_x\}\).
\end{enumerate}
\end{lemma}

\begin{proof}
(1)\( \Rightarrow \)(2): If \( \Delta_{xxx}=0 \)  and  \(  F_{6u} \neq 0\), then, because of Lemmas~\ref{lemu93} and~\ref{lem381}, we have  \( f_{xxx}=a_0 a_x^3+b_0 b_x^3 \) for some \( a_0,b_0\in \C \) and linear forms \( a_x \) and \( b_x \). Since \(  F_{6u}\neq 0 \), \eqref{eqq133} implies that  \( a_0 \) and \( b_0 \) are nonzero, and \( \{a_x,b_x\} \) is linearly independent.

(1)\( \Rightarrow \)(3): If \( \Delta_{xxx}=0 \), then \lref{lemu93} implies that  \( f_{xxx}=a_xb_xc_x \)  for some linearly dependent forms  \( \{a_x,b_x,c_x\}\). Because \(  F_{6u} \neq 0\), \eqref{eq376} implies that each pair of linear forms must be linearly independent.

(2,3)\( \Rightarrow \)(1): If \upitem2 or \upitem3 hold, then \( \Delta_{xxx}=0 \) by \lref{lemu93}, and  \(  F_{6u} \neq 0\) by \eqref{eq376} and \eqref{eqq133}.
\end{proof}

If \( f_{xxx}=a_0 a_x^3+b_0 b_x^3 \) for some \( a_0,b_0\in \C \) and linear forms \( a_x \) and \( b_x \), then, by an easy calculation, \( \theta_{uuxx}=36 a_0 b_0 [abu]^2 a_x b_x \). Assuming that \( \theta_{uuxx} \) is nonzero, this implies that \( a_x \) and \( b_x  \) are uniquely determined by \( f_{xxx} \) (up to order and multiplication by constants). We can use this fact to express \( f_{xxx} \) as a sum of two cubes when that is possible, as seen in the following example.


\begin{example}
Suppose that \( f_{xxx}=2 x_1 (x_1^2 + 6 x_2^2) \). Then \( \Delta_{xxx}=0 \) and so \( f_{xxx} \) is completely reducible. This is no surprise  since  \( f_{xxx}=2 x_1(x_1+i\sqrt6\,x_2)(x_1-i\sqrt6\,x_2) \).

Also \(  F_{6u} = -13824 \,u_3^6\) is nonzero and so, by  \lref{lemq27},  \( f_{xxx}=a_0 a_x^3+b_0 b_x^3 \)  for some linear forms \( a_x \)  and \( b_x \)  such that \( \{a_x,b_x\} \)  is linearly independent, and nonzero \( a_0,b_0\in \C \). To find \( a_x \) and \( b_x \) we factor \( \theta_{uuxx} \):
\[ \theta_{uuxx}=1152\, u_3^2 (x_1^2 - 2 x_2^2)=1152\, u_3^2 (x_1+ \sqrt{2}\, x_2)(x_1- \sqrt{2}\, x_2). \] 
So we can choose \( a_x=x_1+ \sqrt{2}\, x_2 \) and \( b_x=x_1- \sqrt{2}\, x_2 \). Then \( a_0 \) and \( b_0 \) can be determined from \eqref{eqq949} \textup(with \( c_x=x_3 \), for example\textup), or directly, by matching coefficients in \( f_{xxx}=a_0 a_x^3+b_0 b_x^3 \), to get \( a_0=b_0=1 \) and  \[ f_{xxx}=(x_1+ \sqrt{2}\, x_2)^3+(x_1- \sqrt{2}\, x_2 )^3. \] 
\end{example}


%
%
%

\new We see in this example that, not only is \( \theta_{uuxx} \) a reducible quadratic form in \( \xx \), it is also a square as a form in \( \uu \). This is generally true for cubic forms with zero Hessian and leads to other properties of such forms. 


%

\begin{lemma}\label{lemg91}
For a cubic form  \( f_{xxx} \), the following are equivalent:
\begin{enumerate}
\item \( \Delta_{xxx}=0 \)

\item \( \{\partial_1 f_{xxx},\partial_2 f_{xxx},\partial_3 f_{xxx} \} \) is linearly dependent.

\item \( f_{xxz}=0 \) for some nonzero \( (\zz)\in \C^3 \).

\item \( f_{xxx} =g_0 a_x ^3+ h_0 a_x^2b_x+k_0a_xb_x^2+  i_0 b_x^3\) for some linear forms \( a_x \)and \( b_x \), and \( g_0,h_0,k_0,i_0\in \C \).

\item \( f_{xxx}=a_xb_xc_x \) for some linear forms \( a_x \), \( b_x \) and \( c_x \) such that \( \{a_x,b_x,c_x\} \) is linearly dependent.
\end{enumerate}
\end{lemma}

\begin{proof}
\noindent(2)\( \Leftrightarrow \)(3): The equivalence of (2) and (3) is immediate from the expression for \( f_{xxz} \) in \eqref{eqt66}.

\noindent(1)\( \Rightarrow \)(2,3): If  \( f_{xxx} \) is a cube, that is, \( f_{xxx}=a_0 a_x^3 \), then \( f_{xxz}=3 a_0 a_x^2 a_z \), so fixing any nonzero \( (\zz)\in \C^3 \) such that \( a_z= 0 \) we get  \(  f_{xxz}=0 \) for all \( (\xx) \).

Otherwise, by \lref{lem50x},  \( \theta_{uuxx} \) is nonzero and we fix \( (\yy),(\vv)\in \C^3 \ \) so that \( \theta_{vvyy}\in \C \) is  nonzero. Because of the identity 
\[ 4 \theta_{uuyy} \theta_{vvyy} - \theta_{uvyy}^2 =16 \Delta_{yyy}\,\Vert  f_{xxy}\Vert_{x\mapsto[uv]}  \]
and the assumption that \( \Delta_{yyy}=0 \), we have  
\( \theta_{uuyy}= \theta_{uvyy}^2 /(4 \theta_{vvyy}) \), and so  \( \theta_{uuyy} \) is a square. We can write this as \( \theta_{uuyy}= z_0 u_z^2 \) for some fixed nonzero \( (\zz)\in \C^3 \). Moreover, because of the identity
\[  \Vert \theta_{uuyy} \Vert_{u\mapsto\partial f} = 3 f_{xxx} \Delta_{xyy} - f_{xxy} \Delta_{xxy} + f_{xyy} \Delta_{xxx},\] and \( \Delta_{xxx}=0 \), we get \( \Vert u_z ^2\Vert_{u\mapsto\partial f} =0\), which, from \eqref{eq01} and \eqref{eqt66}, is the same as \( f_{xxz}^2=0 \). Thus   \( f_{xxz}=0 \).  

\noindent(3)\( \Rightarrow \)(4): We have  \( f_{xxz}=0 \) for some fixed nonzero \( (\zz)\in \C^3 \). Let \( \{a_x,b_x\} \) be a basis for the two-dimensional vector space \( V=\{u_x\mid u_z=0\} \). Geometrically, \( V \) is the set of lines that pass through the point \( (\zz) \). Since \( \{a_x,b_x\} \) is linearly independent and \( a_z=b_z=0 \), after possible scaling, we can choose this basis so that  \( z_ 1 = a_ 3 b_ 2 - a_ 2 b_ 3\), \(z_ 2 = -a_ 3 b_ 1 + a_ 1 b_ 3\)  and \(
z_ 3 = a_ 2 b_ 1 - a_ 1 b_ 2 \). In this case, \( \J[f_{xxx},a_x,b_x]=f_{xxz}=0 \), so (4) follows from \lref{lemq913}(C1).


\noindent(4)\( \Rightarrow \)(5): Since \( f_{xxx} \) is a binary form in \( a_x \) and \( b_x \), it factors completely. Each of the factors is a linear combination of \( a_x \) and \( b_x \), so the three factors of \( f_{xxx} \) are linearly dependent.

\noindent(5)\( \Rightarrow \)(1): This follows immediately from \eqref{eq370}. 
\end{proof}
The argument in this proof is essentially the same used by   Pasch \cite{pasch1875} in 1875.

\new The most important relationship between a cubic form and its Hessian is the identity
\begin{equation}\label{eqq124}
 \J3[f_{xxx},\Delta_{xxx},u_x^3]=0 
\end{equation}
The left side is a cubic form in \( \uu \)  and has ten coefficients, so  it  can be seen as providing ten linear equations in the ten coefficients of \( \Delta_{xxx} \). The following lemma gives consequences of these equations in a more general context.  

\begin{lemma}\label{lemq115}
Suppose that   \( \J3[f_{xxx},g_{xxx},u_x^3]=0 \) holds for cubic forms  \( f_{xxx} \) and \( g_{xxx} \).
\begin{enumerate}
\item If \( f_{333}\neq 0 \), then \( g_{xxx}=0 \) if and only if \( g_{113}=g_{123}=g_{223}=g_{133}=g_{233}=g_{333}=0 \).

\item If \( f_{333}\neq 0 \) and \(  4 f_{113} f_{223}-f_{123}^2\neq 0 \), then \( g_{xxx}=0 \) if and only if \( g_{113}=g_{123}=g_{223}=g_{333}=0 \).
\end{enumerate}
\end{lemma}

\begin{proof}
Suppose that \( g_{113}=g_{123}=g_{223}=g_{333}=0 \). Setting the coefficients of  \( u_1^3 \), \( u_1^2 u_2 \), \( u_1u_2^2  \) and \( u_2^3 \) to zero in the equation  \( \J3[f_{xxx},g_{xxx},u_x^3]=0\) gives
\begin{align*}
3 f_{333} g_{222} &=- f_{223} g_{233} \\
3 f_{333} g_{122} &=- f_{223} g_{133} - f_{123} g_{233} \\
3 f_{333} g_{112} &=- f_{123} g_{133} - f_{113} g_{233} \\ 
3 f_{333} g_{111} &=- f_{113} g_{133} 
\end{align*}
Since \( f_{333}\neq 0 \), if, in addition, \( g_{133}=g_{233}=0 \), then these equations imply \( g_{111}=g_{112}=g_{122}=g_{222}=0 \) and so \( g_{xxx}=0 \). Thus \upitem1 is proved.

Otherwise, these equations can be used to express \( g_{111} \), \( g_{112} \), \( g_{122} \) and \( g_{222} \) in terms of \( g_{133} \) and \( g_{233} \). Substituting these expressions into \( \J3[f_{xxx},g_{xxx},u_x^3]=0 \) and setting the coefficients of \( u_1 u_3^2 \) and \( u_2 u_3^2 \)  to zero now gives \( ( 4 f_{113} f_{223}-f_{123}^2) g_{233}=0 \) and \( ( 4 f_{113} f_{223}-f_{123}^2) g_{133} =0\). If \(  4 f_{113} f_{223}-f_{123}^2\neq 0 \), this implies \( g_{233}=g_{133}=0 \), and so, with (1), we get  \( g_{xxx}=0 \). 
%
\end{proof}


As an immediate consequence of \eqref{eqq124}  and this lemma we get the following: 
\begin{lemma}
Let \( f_{xxx} \) be a cubic form.
\begin{enumerate}
\item If \( f_{333}\neq 0 \), then \( \Delta_{xxx}=0 \) if and only if \( \Delta_{113}=\Delta_{123}=\Delta_{223}=\Delta_{133}=\Delta_{233}=\Delta_{333}=0 \).

\item If \( f_{333}\neq 0 \) and \(  4 f_{113} f_{223}-f_{123}^2\neq 0 \), then \( \Delta_{xxx}=0 \) if and only if \( \Delta_{113}=\Delta_{123}=\Delta_{223}=\Delta_{333}=0 \).
\end{enumerate}
\end{lemma}
Combined with \lref{lemu93}, we have shown that, generically, four equations, namely,  \( \Delta_{113}=\Delta_{123}=\Delta_{223}=\Delta_{333}=0 \), have to hold so that \( f_{xxx}=a_xb_xc_x \) for some linear forms \( a_x \), \( b_x \) and \( c_x \) such that \([abc]=0\). 

\section{Complete Reducibility}\label{completereducibility}

Suppose for the moment that \(f_{xxx} \) is completely reducible, but we don't know its factorization. That is,  \(f_{xxx}= a_xb_xc_x \), but we don't know \(a_x \), \(b_x\) or \(c_x\).  In this circumstance we can use the identity \( f_{xyy}= a_y b_y c_x+a_y  c_yb_x +b_yc_ya_x \) to find the factors of \(f_{xxx}\) as follows.

Suppose that \((\yy)\in \C^3\) is a zero of \( f_{xxx}\),  that is,  \( f_{yyy}=a_yb_yc_y=0\). Then, except in special cases, we have \( a_y=0\), \(b_y\neq 0\) and \(c_y \neq 0\), or  \( a_y\neq 0\), \(b_y= 0\) and \(c_y \neq 0\), or  \( a_y\neq 0 \), \(b_y\neq 0\) and \(c_y =0\).  Hence \(  f_{xyy}\) is \( b_y c_y a_x\), \( a_y  c_y b_x\) or \(a_y b_y c_x\). Either way, \(f_{xyy}\) is a nonzero linear factor of \(f_{xxx}\).  A geometer would recognize that  \( (\yy) \) is a point on the curve \( f_{xxx}=0 \) and \( f_{xyy}=0 \) is the equation of the tangent line to the curve at that point.  To find all three of the linear factors of \( f_{xxx}\), we need to find three zeros of \(f_{xxx}\) and the three corresponding tangent lines. 

Let's look for three zeros of \(f_{xxx}\) along the line that joins two fixed points \( (\yy),(\zz)  \in \C^3 \). The equation of the line passing through those points is simply \( [xyz]=0 \), or \( u_x=0 \) where the (fixed) line coordinates \( \uu \) are determined by the substitution \( u\mapsto [yz] \) defined in  \eqref{equ92}. Points on the line are parametrized by the substitution  \( x\mapsto X \) defined by 
\begin{equation}\label{eq711}
x_1\mapsto X_1 y_1+X_2 z_1 \quad x_2\mapsto X_1 y_2+X_2 z_2 \quad x_3\mapsto X_1 y_3+X_2 z_3 ,
\end{equation}
where \( X_1 \) and \( X_2 \) are new variables. The values of the cubic form \( f_{xxx} \) at points on the line are given by  (compare \eqref{eqs11})
\begin{equation}\label{eqw97}
f_{XXX}=\Vert f_{xxx}\Vert_{x\mapsto X}=f_{yyy} X_1^3+f_{yyz} X_1^2 X_2^{\phantom{2}}+f_{yzz} X_1^{\phantom{2}} X_2^2+ f_{zzz} X_2^3
\end{equation}
Being a binary form, \( f_{XXX} \) factors completely, that is, there are \( A_1 \), \( A_2 \), \( B_1 \), \( B_2 \), \( C_1 \), \( C_2\in \C \) such that
\begin{equation}\label{eq75}
f_{XXX}=(A_1 X_1 + A_2 X_2 ) (B_1 X_1 + B_2 X_2 ) (C_1 X_1 + C_2 X_2 ).
\end{equation}
One of the zeros of \( f_{XXX} \) occurs when \( X_1=A_2 \), \( X_2=-A_1 \). This corresponds to the point \( ( A_2y_1-A_1 z_1,A_2 y_2-A_1 z_2, A_2 y_3-A_1 z_3) \) which is therefore a zero of \( f_{xxx} \), and geometrically, an intersection point of the curve \( f_{xxx}=0 \) and the line \( u_x =0 \). Of course, the other factors of \( f_{XXX} \) give other intersection points. The equations of the tangent lines to the curve at these points are obtained from  \( f_{xyy} \) by setting \( (\yy) \) to the coordinates of the intersection  points: 
%
\begin{equation}\label{eq73}
\begin{aligned}
 A_x &= f_{xyy} A_2^2 - f_{xyz} A_1 A_2 + f_{xzz} A_1^2\\
 B_x &= f_{xyy} B_2^2 - f_{xyz} B_1 B_2 + f_{xzz} B_1^2\\
C_x &= f_{xyy} C_2^2 - f_{xyz} C_1 C_2 + f_{xzz} C_1^2
\end{aligned}
\end{equation}
Now set \( L_{xxx}=A_xB_xC_x \), an obviously completely reducible cubic form. By matching coefficients in \eqref{eqw97} and \eqref{eq75} we get
\begin{align*}
f_{yyy} &= A_1 B_1 C_1 & f_{yyz} &= A_2 B_1 C_1 + A_1 B_2 C_1 + A_1 B_1 C_2\\
 f_{zzz} &= A_2 B_2 C_2 & f_{yzz} &= A_2 B_2 C_1 + A_2 B_1 C_2 + A_1 B_2 C_2 
\end{align*}
Using these equations,  \( L_{xxx} \) can be written entirely in terms of the polars of \( f_{xxx} \):
\begin{equation}\label{eqg845}
\begin{aligned}
L_{xxx} &= f_{xzz}^3 f_{yyy}^2 - f_{xyz} f_{xzz}^2 f_{yyy} f_{yyz} + f_{xyy} f_{xzz}^2 f_{yyz}^2 + f_{xyz}^2 f_{xzz} f_{yyy} f_{yzz} \\
&\qquad - 2 f_{xyy} f_{xzz}^2 f_{yyy} f_{yzz} - f_{xyy} f_{xyz} f_{xzz} f_{yyz} f_{yzz} +   f_{xyy}^2 f_{xzz} f_{yzz}^2  \\
&\qquad - f_{xyz}^3 f_{yyy} f_{zzz} +   3 f_{xyy} f_{xyz} f_{xzz} f_{yyy} f_{zzz} + f_{xyy} f_{xyz}^2 f_{yyz} f_{zzz}  \\
&\qquad -  2 f_{xyy}^2 f_{xzz} f_{yyz} f_{zzz} - f_{xyy}^2 f_{xyz} f_{yzz} f_{zzz} + f_{xyy}^3 f_{zzz}^2
\end{aligned}
\end{equation}

We might expect that, if the discriminant of \(f_{XXX}\), 
\[ \delta=f_{yyz}^2 f_{yzz}^2 - 4 f_{yyy} f_{yzz}^3 - 4 f_{yyz}^3 f_{zzz} + 18 f_{yyy} f_{yyz} f_{yzz} f_{zzz} - 27 f_{yyy}^2 f_{zzz}^2 \]
  is nonzero so that the linear factors of  \(f_{XXX}\) are distinct, then \(A_x\), \(B_x\), and \(C_x\) are distinct linear factors of \(f_{xxx}\), and  \(f_{xxx}\) is just a scalar multiple of \(L_{xxx}\). This is made clear by the fact that 
 \begin{equation}\label{eqg847}
 L_{xxx}+\delta f_{xxx}=0
 \end{equation}

\new What if \(f_{xxx} \) is not completely reducible? Given an arbitrary cubic form \(f_{xxx}\)  we can still choose points \((\yy),(\zz)\in \C^3\) and construct \(L_{xxx}\) as above.  Geometrically, \( L_{xxx} \) is the product of the tangent lines to the three intersection points of the curve \( f_{xxx}=0 \) and the line \( u_x=0 \).  By construction \(L_{xxx}\) is a completely reducible form, even if \(f_{xxx}\) is not.  The equations \eqref{eqw97}-\eqref{eqg845} above are still true, but, in the general situation, \eqref{eqg847} becomes the identity
\begin{equation}\label{eqq127}
L_{xxx}+\delta f_{xxx}=[xyz]^2 \Vert \Gamma_{4ux} \Vert_{u\mapsto[yz]}
\end{equation}
where  \( \Gamma_{4ux} \) is defined in \eqref{eq837}. From this equation we see immediately that, if  \((\yy),(\zz)\in \C^3\) are chosen so that \( \delta \) is nonzero, and \( \Vert \Gamma_{4ux} \Vert_{u\mapsto[yz]}=0 \), then \( f_{xxx} \) is a multiple of \( L_{xxx} \) and hence completely reducible.
 Conveniently,  the discriminant of \( f_{XXX} \) is related to \( F_{6u} \) by the identity
 \begin{equation}\label{eq54}
 \delta = \Vert  F_{6u}\Vert_{u\mapsto[yz]}
 \end{equation}
This equation is an example of the Clebsch Transfer Principle \cite{gundelfinger1873}, \cite[Article~215]{grace} which relates concomitants of binary forms, such as \( \delta \),  to certain concomitants of ternary forms,  such as \( F_{6u} \).  


Because of \eqref{eq54}, \( \delta \) can be replaced by \(  F_{{6u}} \) in \eqref{eqq127}:
\begin{equation}\label{eqq10}
L_{xxx}+ \Vert  F_{6u} \Vert_{u\mapsto[yz]}f_{xxx}=  [xyz]^2 \Vert \Gamma_{4ux} \Vert_{u\mapsto[yz]}
\end{equation}

\begin{lemma}\label{lemq21}
For a cubic form \( f_{xxx} \), suppose that \( \uu\in \C \) are fixed so that \(  F_{6u} \) is nonzero. Then the following are equivalent:
\begin{enumerate}
\item \( f_{xxx} \) is completely reducible.

\item \( \Delta_{xxx} \) is a multiple of \( f_{xxx} \). 

\item \( \Pi_{u4x}=0 \)

\item \( \Gamma_{4ux}=0 \)
\end{enumerate}
\end{lemma}

\begin{proof}
If  \( f_{xxx} \) is completely reducible, then, by \eqref{eq370}, \( \Delta_{xxx} \) is a multiple of \( f_{xxx} \), so the antisymmetry of the Jacobian ensures that  \( \J[f_{xxx},\Delta_{xxx},u_x]=0 \). Hence, straight from the definition \eqref{eq837}, we get  \( \Pi_{u4x}=0 \).  

If \( \Pi_{u4x}=0 \), then \( \Gamma_{4ux}=0 \) follows directly from the definition \eqref{eq837}.

Finally, suppose that \( \Gamma_{4ux}=0 \). Choose \( \yy\in \C \) and \( \zz\in\C \) so that \( u_x=[xyz] \). Then the value of the form \(  F_{6u} \) at the fixed  \( \uu\in \C \) is equal to the value of the form \( \Vert  F_{6u}\Vert_{u\mapsto[yz]}\) at the fixed \( \yy, \zz\in\C \), in particular, this value is nonzero.  Since \( \Gamma_{4ux}=0 \),  \eqref{eqq10} becomes
\( F_{6u} f_{xxx}=-L_{xxx}\) which, because \( L_{xxx} \) is completely reducible, implies \( f_{xxx} \) is completely reducible.
\end{proof}

It is remarkable that in a footnote in Gundelfinger's 1871 paper \cite[p.~227]{gundelfinger1871c} he claims, in effect,  that  it is easy to show that \(  F_{6u} f_{xxx} - u_x^2  \Gamma_{4ux} \) is the product of the tangent lines to the three intersection points of the curve \( f_{xxx}=0 \) and the line \( u_x=0 \).
On the following page, he shows that if \( \Delta_{xxx} \) is a multiple of \( f_{xxx} \), then \( \Gamma_{4ux} =0 \). So he was very close to proving \lref{lemq21} as we have done above. Unfortunately, in the rest of the paper, he  took a more convoluted path in his attempt to prove the same result.

\new Except in very special cases (\lref{lem381}),  for a given  \( f_{xxx} \),  \(  F_{6u} \neq 0 \) will hold for almost all choices of  \( \uu\in \C \), and once \( \uu \) are chosen satisfying this condition, \( f_{xxx} \) is completely reducible if and only if the three coefficients of  \( \xx  \)  in \( \Gamma_{4ux}  \) are zero. 
The fact that, generically, three conditions should suffice to determine whether a cubic form is completely reducible was well-known in the 19th century \cite[\S2]{aronhold1849}, \cite{thaer1879}. 

For comparison, if  \( \uu\in \C \) are fixed, then \( \Pi_{u4x} \) is a degree \( 4 \) form in \( \xx \), so to confirm that \( \Pi_{u4x} =0 \),  \( 15 \) coefficients have to be checked.



\new In \lref{lemq21},  the variables  \( \uu \) are fixed such that \( F_{6u} \) is nonzero. This restriction on \( \uu \) is needed only to prove that (4) implies (1) and can be easily removed so long as we require the conditions in  \lref{lemq21} to hold for all \( \uu \), as in the following theorem.

\begin{theorem}\label{thmq111}
For a cubic form \( f_{xxx} \), the following are equivalent:
\begin{enumerate}
\item \( f_{xxx} \) is completely reducible.

\item \( \Delta_{xxx} \) is a multiple of \( f_{xxx} \). 

\item \( \Pi_{u4x}=0 \)

\item \( \Gamma_{4ux}=0 \)
\end{enumerate}
\end{theorem}

\begin{proof}
The only thing to add to the proof of \lref{lemq21} is the following: Suppose that  \( \Gamma_{4ux}=0 \). There are two cases. If \(  F_{6u}=0 \), then \( f_{xxx} \) is completely reducible because of \lref{lem381}. Otherwise, we can fix \( \uu\in \C \) so that \(  F_{6u}\in \C \) is nonzero. Then  \( \Gamma_{4ux} \) is zero for those \( \uu\in \C \) and so \( f_{xxx} \) is completely reducible by \lref{lemq21}.
\end{proof}

The advantage of \tref{thmq111} over \lref{lemq21} is that (1) and (2) are equivalent conditions, independent of the values of \( F_{6u} \). The disadvantage is that  \( \Pi_{u4x}\) and  \( \Gamma_{4ux}=0 \) are treated as forms in both \( \xx \) and \( \uu \), with \( 45  \) coefficients each. So to confirm (3) and (4), all \( 45 \) coefficients have to be checked. 

In fact, the coefficients of \( \Pi_{u4x}\) and  \( \Gamma_{4ux}=0 \) are not linearly independent. The identities \( \con_{ux}[ \Pi_{u4x}]=0 \) and \( \con_{ux}[\Gamma_{4ux}]=0 \) give 10 linear equations satisfied by the coefficients of each form. So it suffices to check only \( 35 \) coefficients to confirm that \( \Pi_{u4x}=0 \) or  \( \Gamma_{4ux}=0 \).

\new It seems that the only result from previous sections needed  to prove  \lref{lemq21} and \tref{thmq111}  is just one part of \lref{lem381}: If \( F_{6u}=0 \), then \( f_{xxx} \) is completely reducible. Unfortunately, there seems to be no direct way of proving this that avoids the  discussion of the conditions \( \theta_{uuxx}=0 \) and \( \Delta_{xxx}=0 \) found in Sections~\ref{cubic} and~\ref{hessian}.

\new It is well worth noticing that our discussion provides a method for factoring \( f_{xxx} \) when it is completely reducible and \( F_{6u}\neq 0 \): 


\bigskip \textit{Fix \( (\uu)\in \C^3 \) so that \( F_{6u}\in \C \) is nonzero. Fix \( (\yy)  \) and \( (\zz) \) so that 
\( u_1=  y_2 z_3-y_3 z_2  \),    \( u_2= y_3 z_1 - y_1 z_3 \) and    \( u_3 =  y_1 z_2-y_2 z_1 \).
Calculate \( f_{XXX} \) from \eqref{eqw97} and factor it as in \eqref{eq75}. Calculate \( A_x \), \( B_x \) and \( C_x \) from  \eqref{eq73}. Then, from \eqref{eqg847} and \eqref{eq54}, we get  \[  f_{xxx}=- \dfrac1{F_{6u} }A_xB_xC_x. \]}

\textit{Note that merely choosing \( (\yy) \) and \( (\zz) \) to be arbitrary distinct points on the line \( u_x=0 \) would give same linear factors of \( f_{xxx} \), just multiplied by some scalar. }

\new Let's work through an example.

\begin{example}
Let
\[ f_{xxx}=x_1^3 - 6 x_1 x_2^2 - 6 x_2^3 + 6 x_1^2 x_3 + 18 x_1 x_2 x_3 + 12 x_2^2 x_3 + 
 4 x_3^3.  \] Since \( \Delta_{xxx}= -108 f_{xxx} \), this form is completely reducible by \tref{thmq111}. \( F_{6u} \) is rather complicated, but it is easy to choose \( (\uu) \) so that \( F_{6u} \) is nonzero. For example, if we fix \( (\uu)=(0,0,1) \), then \( F_{6u}=-108 \). Now choose \( (\yy)=(1,0,0) \) and \( (\zz)=(0,1,0) \). From \eqref{eqw97} we get   \( f_{XXX}=X_1^3 - 6 X_1 X_2^2 - 6 X_2^3 \). Let \( \alpha \), \( \beta \) and \( \gamma \) be the roots of \( x^3-6x-6 \). Then \( f_{XXX} \) factors as in \eqref{eq75} with \( A_1=B_1=C_1=1 \), \( A_2=-\alpha \), \( B_2=-\beta \) and \( C_2=-\gamma \). The linear factors of \( f_{xxx} \) given by \eqref{eq73} are 
 \begin{align*}
 A_x&=3 ( \alpha^2-2 ) x_1 - 6 (3 + 2 \alpha) x_2 + 6 (1 + \alpha) (2 + \alpha) x_3\\
 B_x&=3 ( \beta^2-2) x_1 - 6 (3 + 2 \beta) x_2 + 6 (1 + \beta) (2 + \beta) x_3\\
 C_x&=3 ( \gamma^2-2 ) x_1 - 6 (3 + 2 \gamma) x_2 + 6 (1 + \gamma) (2 + \gamma) x_3
 \end{align*}
 and finally, \( f_{xxx}= \dfrac1{108} A_xB_xC_x \).
 
Starting all over again, a ``better'' choice might be \( (\uu)=(1,1,0) \) with \( (\yy)=(-1,1,0) \) and \( (\zz)=(0,0,1) \) because  then \( F_{6u}=-432 \) is nonzero and \( f_{XXX}=-X_1^3 + 4 X_2^3 \) is easy to factor with \( A_1=B_1=C_1=-1 \), \( A_2=\sqrt[3]4 \), \( B_2=\omega\sqrt[3]4  \) and \( C_2=\omega^2\sqrt[3]4  \). The linear factors of \( f_{xxx} \) given by \eqref{eq73}  are now
 \begin{align*}
 A_x&=6 (\sqrt[3]4 - \sqrt[3]2) x_1 + 6 (\sqrt[3]4 - 2 \sqrt[3]2) x_2 + 12 x_3\\
 B_x&=6 (\omega \sqrt[3]4 - \omega^2 \sqrt[3]2) x_1 + 6 (\omega \sqrt[3]4 - 2 \omega^2 \sqrt[3]2) x_2 + 12 x_3\\
 C_x&=6 (\omega^2 \sqrt[3]4 - \omega \sqrt[3]2) x_1 + 6 (\omega^2  \sqrt[3]4 - 2 \omega \sqrt[3]2) x_2 +  12 x_3
 \end{align*}
and \( f_{xxx}=\dfrac1{432}A_xB_xC_x \).

Because of the unique factorization property of polynomials, the two factorizations of \( f_{xxx} \) we have just found should be the same up to multiplication of constants. The first step in confirming that would to be to calculate the roots of  \( x^3-6x-6 \) using Cardano's formula to get  \( \alpha=\sqrt[3]2+\sqrt[3]4 \), \( \beta=\omega\sqrt[3]2+\omega^2\sqrt[3]4 \) and \( \gamma=\omega^2\sqrt[3]2+\omega\sqrt[3]4 \).
\end{example}

\new Here are two other necessary and sufficient conditions for a cubic form \( f_{xxx} \) to be completely reducible.

\begin{lemma}\label{lemq120}
For a cubic form \( f_{xxx} \), the following are equivalent:
\begin{enumerate}
\item \( f_{xxx} \) is completely reducible.

\item \( \Delta_{xxx}^2-S f_{xxx}^2=0 \)

\item \( F_{6u} \Delta_{xxx} - S_{uuu}^2 f_{xxx }=0\)
\end{enumerate}
\end{lemma}

\begin{proof}
If \( f_{xxx} \) is completely reducible, then (2) and (3) hold because of equations \eqref{eq370}-\eqref{eq376}.

If (2) holds then \(  (\Delta_{xxx}+\sqrt{S} f_{xxx})(\Delta_{xxx}-\sqrt{S} f_{xxx})=0 \) and so either  \( \Delta_{xxx}= \sqrt{S} f_{xxx} \) or  \( \Delta_{xxx}=- \sqrt{S} f_{xxx} \). Either way, \(  \Delta_{xxx} \) is a multiple of \( f_{xxx} \) and  \( f_{xxx} \) is completely reducible by \tref{thmq111}.

Suppose that (3) holds.  If \( F_{6u}=0 \), then \( f_{xxx} \) is completely reducible by \lref{lem381}. Otherwise, fix \( \uu\in \C \) so that \( F_{6u}\in \C \) is nonzero. Then \( F_{6u} \Delta_{xxx} - S_{uuu}^2 f_{xxx }=0\) shows that \( \Delta_{xxx} \) is a multiple of \( f_{xxx} \), and  \( f_{xxx} \) is completely reducible  by \tref{thmq111}. 
\end{proof}

As explained earlier, some of the conditions in \eqref{eqy5} are necessary and sufficient for a cubic form \( f_{xxx} \) to be completely reducible so long as \( S\neq 0 \).

\begin{lemma}\label{lemq121}
For a cubic form \( f_{xxx} \) with \( S\neq 0 \), the following are equivalent:
\begin{enumerate}
\item \( f_{xxx} \) is completely reducible.

\item \( S \Delta_{xxx}-Tf_{xxx}=0 \)

\item   ﻿﻿\(   T \Delta_{xxx}-S^2 f_{xxx} =0  \)

\item \( S F_{6u} -  S_{uuu} T_{uuu}=0 \)

\item \( \J2[\Delta_{xxx},  \Delta_{xxx}, u_x^2] - 4 S \theta_{uuxx}=0 \)
\end{enumerate}
\end{lemma}

\begin{proof}
If \( f_{xxx} \) is completely reducible, then (2)-(5) hold because of equations \eqref{eq370}-\eqref{eq376}. If (2) or (3) hold, then \( \Delta_{xxx} \) is a multiple of \( f_{xxx} \) and so \( f_{xxx} \) is completely reducible by \tref{thmq111}. If (4) holds with \( S\neq 0 \), then the identity 
\[ 432 \,S\, \Gamma_{4ux} = 6 \con_{ux}^2[(S F_{6u} -  S_{uuu} T_{uuu}) f_{xxx}] - 
 \con_{ux}^3[(S F_{6u} -  S_{uuu} T_{uuu}) f_{xxx}] u_x \] 
 implies that \( \Gamma_{4ux}=0 \), so \( f_{xxx} \) is completely reducible  by \tref{thmq111}. If (5) holds, then the identity
 \[ S \Delta_{xxx}-Tf_{xxx}=-96 \con_{ux}^2[(\J2[\Delta_{xxx},  \Delta_{xxx}, u_x^2] - 4 S \theta_{uuxx})f_{xxx}] \] implies that (2) holds and hence \( f_{xxx} \) is completely reducible.
 \end{proof}

\section{Special Case I}


In this section we consider the special case of \lref{lemq21} with  \( (\uu)=(0,0,1) \), that is, \( u_x=x_3 \).  In this circumstance, 
\[  F_{6u}= f_{112}^2 f_{122}^2 - 4 f_{111} f_{122}^3 - 4 f_{112}^3 f_{222} + 18 f_{111} f_{112} f_{122} f_{222} - 27 f_{111}^2 f_{222}^2 \]
is the discriminant of \( f_{111} x_1^3 + f_{112} x_1^2 x_2 + f_{122} x_1 x_2^2 + f_{222} x_2^3 \), and 
\begin{align*}
288\, \Gamma_{4ux}  &=6 \left(\Delta_{122} \theta_{3311} - \Delta_{112} \theta_{3312} + 3 \Delta_{111} \theta_{3322}\right)x_1 \\
&\qquad + 
 6 \left(3 \Delta_{222} \theta_{3311} - \Delta_{122} \theta_{3312} + \Delta_{112} \theta_{3322}\right) x_2 \\[2pt]
 &\qquad -   (\Delta_{122} \theta_{1311} - \Delta_{112} \theta_{1312} + 3 \Delta_{111} \theta_{1322} + 3 \Delta_{222} \theta_{2311}  - \Delta_{122} \theta_{2312} \\[2pt]
 &\qquad\quad+ \Delta_{112} \theta_{2322} - 4 \Delta_{223} \theta_{3311} + 2 \Delta_{123} \theta_{3312} - 4 \Delta_{113} \theta_{3322})\, x_3
\end{align*}
For compactness we are writing  \( \Gamma_{4ux} \) in terms of the coefficients of \( \Delta_{xxx} \) and \( \theta_{uuxx} \) using the identity
\[ 2304\, \Gamma_{4ux} =14 \J2[\Delta_{xxx}, \theta_{uuxx}, u_x^2] -  \con_{ux}[\J2[\Delta_{xxx}, \theta_{uuxx}, u_x^2] u_x]. \]
Right from the definition \eqref{eq837}, \( F_{6u} \) can also be expressed in terms of the coefficients of \( \theta_{uuxx} \),
\[ 48 F_{6u}= 4 \theta_{3311} \theta_{3322}-\theta_{3312}^2, \] 
so the  \( (\uu)=(0,0,1) \) case of  \lref{lemq21} can be expressed as follows.

\begin{lemma}
Let \( f_{xxx} \) be a ternary cubic form such that
\( 4 \theta_{3311} \theta_{3322}-\theta_{3312}^2 \) is nonzero.
Then \( f_{xxx} \) is completely reducible if and only if 
\begin{align*}
 \Delta_{122} \theta_{3311} - \Delta_{112} \theta_{3312} + 3 \Delta_{111} \theta_{3322}&=0 \\
3 \Delta_{222} \theta_{3311} - \Delta_{122} \theta_{3312} + \Delta_{112} \theta_{3322}&=0
\end{align*}
\vspace{-20pt}
\begin{multline*}
\Delta_{122} \theta_{1311} - \Delta_{112} \theta_{1312} + 3 \Delta_{111} \theta_{1322}  
+ 3 \Delta_{222} \theta_{2311}  - \Delta_{122} \theta_{2312}\\
+ \Delta_{112} \theta_{2322} 
- 4 \Delta_{223} \theta_{3311} + 2 \Delta_{123} \theta_{3312} - 4 \Delta_{113} \theta_{3322}=0
\end{multline*}

\end{lemma}

In 1910, Glenn \cite[p.~89]{glenn} proved this lemma (with the unneeded extra assumption that \( f_{111}\neq 0 \))  using completely different methods.  In 1930, Copeland \cite{copeland} derived necessary and sufficient conditions for the complete reducibility of a ternary cubic form in this same special case.  Copeland's  reducibility condition is expressed in terms of the rank of a certain \( 7\times 9 \) matrix with entries coming from the coefficients of the form.  Specifically, the form is completely reducible if and only if the rank of the matrix is less than seven.

\new A much simpler case occurs when  \( f_{111}=f_{222}=0 \) and   \( f_{112} f_{122}\neq 0 \) since
\begin{align*}
 F_{6u} &=f_{112}^2f_{122}^2\neq 0 \\[2pt]
\Gamma_{4ux} &=f_{122}^2 (f_{113}^2 f_{122} - f_{112} f_{113} f_{123} + f_{112}^2 f_{133})\,x_1 \\
&\quad+ f_{112}^2 ( f_{112} f_{223}^2-f_{122} f_{123} f_{223}  + f_{122}^2 f_{233}) \,x_2 \\
&\quad+ (f_{113}^2 f_{122}^2 f_{223} - f_{112} f_{113} f_{122} f_{123} f_{223} + f_{112}^2 f_{113} f_{223}^2 + f_{112}^2 f_{122}^2 f_{333})\, x_3\\[2pt]
L_{xxx}&=-(f_{112 }x_2 + f_{113} x_3) (f_{122} x_1 + 
   f_{223} x_3) \\
   &\qquad\cdot\left(f_{112} f_{122} (f_{112} x_1 + f_{122} x_2) + ( 
      f_{112} f_{122} f_{123} - f_{112}^2 f_{223}-f_{113} f_{122}^2 ) x_3\right)
\end{align*}

\begin{corollary}
Let \( f_{xxx} \) be a ternary cubic form  with \( f_{111}=f_{222}=0 \) and   \( f_{112} f_{122}\neq 0 \). 
Then \( f_{xxx} \) is completely reducible if and only if 
\begin{gather*}
 f_{113}^2 f_{122} - f_{112} f_{113} f_{123} + f_{112}^2 f_{133}  =0\qquad
 f_{112} f_{223}^2-f_{122} f_{123} f_{223}  + f_{122}^2 f_{233} =0\\
f_{113}^2 f_{122}^2 f_{223} - f_{112} f_{113} f_{122} f_{123} f_{223} + f_{112}^2 f_{113} f_{223}^2 + f_{112}^2 f_{122}^2 f_{333} =0
\end{gather*}
When these equations hold,
\begin{align*}
 f_{xxx}&=\frac1{f_{112}^2f_{122}^2}(f_{112 }x_2 + f_{113} x_3) (f_{122} x_1 +   f_{223} x_3) \\
 &\qquad\cdot (f_{112} f_{122} (f_{112} x_1 + f_{122} x_2) + ( 
      f_{112} f_{122} f_{123} - f_{112}^2 f_{223}-f_{113} f_{122}^2 ) x_3)
\end{align*} 
\end{corollary}

\section{Special Case II}

We consider the complete reducibility of \( f_{xxx} \) under the condition that \( f_{333}\neq 0 \).  By  \tref{thmq111}, \( f_{xxx} \) is completely reducible, if and only if \( \Delta_{xxx} \) is a multiple of \( f_{xxx} \), if and only if  the cubic form \( g_{xxx}=\Delta_{333}f_{xxx}-f_{333}\Delta_{xxx} \) is zero. The form \( g_{xxx} \) has ten coefficients, but, since \( g_{333}=0 \)  by construction, it seems that only nine have to be checked to see if \( f_{xxx} \) is completely reducible. In fact we will show that it suffices to check only five coefficients of \( g_{xxx} \), and, with an additional weak assumption, only three. 

We have already noted the important identity  \( \J3[f_{xxx},\Delta_{xxx},u_x^3] =0 \). Since, in addition,  \( \J3[f_{xxx},f_{xxx},u_x^3]=0\), we have \( \J3[f_{xxx},g_{xxx},u_x^3]=0 \).   By construction, \( g_{333}=0 \), so directly from \lref{lemq115}, we have new criteria for the complete reducibility of \( f_{xxx} \). 

\begin{lemma}\label{lem113}
Let \( f_{xxx} \) be a ternary cubic form such that \( f_{333}\neq 0 \). Set  \[ g_{xxx}=\Delta_{333}f_{xxx}-f_{333}\Delta_{xxx} .\]

\begin{enumerate}
\item  \( f_{xxx} \) is completely reducible if and only if \( g_{113}=g_{123}=g_{223}=g_{133}=g_{233}=0 \). 

\item If, in addition,  \( 4 f_{113} f_{223}-f_{123}^2 \neq 0 \), then  \( f_{xxx} \) is completely reducible if and only if  \( g_{113}=g_{123}=g_{223}=0 \). 
\end{enumerate}
\end{lemma}

The claim that the  five equations \( g_{113}=g_{123}=g_{223}=g_{133}=g_{233}=0 \)  are necessary and sufficient for the complete reducibility of \( f_{xxx} \) was  proved by Brill in 1893 \cite{brill1893} and also by Junker \cite{junker} in 1894 using independent methods.

\begin{example}
Suppose that  \( f_{xxx}=x_3^3 + x_3 x_1^2 + x_2^3 \). Then \( f_{333}=1\neq 0  \), \( \Delta_{333}=0 \) and \( g_{xxx}=\Delta_{xxx}= -12 x_1^2 x_2 + 36 x_2 x_3^2\).   Since \( g_{xxx} \) is not zero, \( f_{xxx} \) is not completely reducible. In fact, by \lref{lem113}\upitem1, to know that \( f_{xxx} \) is not completely reducible it suffices to see that one the five coefficients, \( g_{113}\), \( g_{123}\), \(g_{223}\), \(g_{133}\), \(g_{233} \), is nonzero, namely \( g_{233}=-36 \). 

Note that \lref{lem113}\upitem2 does not apply  because  \( 4 f_{113} f_{223}-f_{123}^2 = 0 \), and, even though \( g_{113}=g_{123}=g_{223}=0 \), \( f_{xxx} \) is not completely reducible. 
\end{example}

One further special case investigated by Brioschi \cite{brioschi1876} in 1876 is \( f_{333}\neq 0 \) and \( f_{133}=f_{233}=0 \).  In this situation, \( \Delta_{333}=3 ( 4 f_{113} f_{223}-f_{123}^2 ) f_{333} \) and so  \[ g_{xxx}=f_{333}\left(3 ( 4 f_{113} f_{223}-f_{123}^2 ) f_{xxx} - \Delta_{xxx}\right) \] Since \( f_{333} \) is nonzero, the equations \( g_{113}=g_{123}=g_{223}=0 \) simplify to three equations that are degree three in the coefficients of \( f_{xxx} \):
\begin{equation}\label{eqq119}
\begin{aligned}
f_{113} ( 4 f_{113} f_{223}-f_{123}^2) + 3 (f_{112}^2 - 3 f_{111} f_{122}) f_{333}&=0\\
f_{123} ( 4 f_{113} f_{223}-f_{123}^2 ) + 3 (f_{112} f_{122} - 9 f_{111} f_{222}) f_{333}&=0\\
f_{223} (4 f_{113} f_{223}-f_{123}^2 ) + 3 (f_{122}^2 - 3 f_{112} f_{222}) f_{333}&=0
\end{aligned}
\end{equation}
\lref{lem113}(2) now can be expressed as the following.

\begin{lemma}
Suppose that  \( f_{333}\neq 0 \), \( \Delta_{333}\neq 0 \) and \( f_{133}=f_{233}=0 \). Then \( f_{xxx} \) is completely reducible if and only if \eqref{eqq119} holds.
\end{lemma}

It seems that  Brioshi's case is rather special so it is worth noticing that, if \( f_{xxx} \) is an arbitrary cubic form with \( f_{333}\neq 0 \), then the linear substitution
\[ x_1\mapsto x_1\qquad x_2\mapsto x_2\qquad x_3 \mapsto x_3 - \frac{f_{133}}{3 f_{333}}  x_1 - \frac{f_{233}}{3 f_{333}}x_2 \] gives a cubic form for which the coefficient of \( x_3^3 \) is unchanged, and the coefficients of \( x_1x_3^2 \) and \( x_2x_3^3 \) are zero.
Hence this discussion can be applied to the transformed cubic form.

\section{Cubic Forms with Symmetry}

In this section, we investigate the complete reducibility of cubic forms \( f_{xxx} \) that are unchanged by even permutations of \( \xx \). Such forms are interesting because there is a simple test for their complete reducibility and we can carry out the factorization when it occurs. 



Before discussing the general case, we note one frequently occurring special case: \( f_{xxx}= x_1^3+x_2^3+x_3^3-3 x_1 x_2 x_3\). Since \( \Delta_{xxx}=-27 f_{xxx} \) this form is completely reducible by \tref{thmq111}. Indeed
\begin{equation}\label{eqq117}
f_{xxx}=(x_1+x_2+x_3)(x_1+\omega  x_2+ \omega^2 x_3)( x_1+\omega^2 x_2+\omega x_3)
\end{equation}
where \( \omega=e^{2\pi i /3} \).

\new For the general case,  suppose that \( f_{xxx} \) is a cubic form that is unchanged by even permutations of \( \xx \).  It is not hard to show that \( f_{xxx} \) can be written as
\begin{equation}\label{eqq113}
\begin{aligned}
 f_{xxx}&=  a(2x_1-x_2-x_3)(2x_2-x_3-x_1)(2x_3-x_1-x_2)\\
 &\qquad+ b(x_1-x_2)(x_2-x_3)(x_3-x_1)\\
 &\qquad+c (x_1+x_2+x_3)^3+ d (x_1^3+x_2^3+x_3^3-3 x_1x_2x_3)
\end{aligned}
\end{equation}
for some \( a ,b,c,d\in \C \). Writing \( f_{xxx} \) this way, as a linear combination of completely reducible forms, simplifies the upcoming algebra. 

The Hessian of \( f_{xxx} \) is 
\[ \Delta_{xxx}=81 d^2 f_{xxx}-108\left((27a^2+b^2)c+d^3\right)(x_1^3+x_2^3+x_3^3-3 x_1x_2x_3)  \] 
By \tref{thmq111}, \( f_{xxx} \)  is completely reducible if and only if \( \Delta_{xxx} \) is a multiple of \( f_{xxx} \). This can happen in two ways: Either \( f_{xxx} \) is a multiple of \( x_1^3+x_2^3+x_3^3-3 x_1x_2x_3 \), or \( (27a^2+b^2)c+d^3=0 \). In the  first case \( f_{xxx} \) factors as in \eqref{eqq117},   so it remains to consider  the factorization of \( f_{xxx} \) in the case \( (27a^2+b^2)c+d^3=0 \).

%

\new We will  follow the argument in Section~\ref{completereducibility} to find the factors of \( f_{xxx} \). First we notice that the last two terms of \( f_{xxx} \) are zero if \( x_1+x_2+x_3=0 \). Consequently,  the concomitants of \( f_{xxx}  \), which are very complicated when written in full, become much simpler when \( (\uu)=(1,1,1)  \). For example, with that choice for \( \uu \), we get \(  F_{6u}= 729 (27 a^2 + b^2)^2\). 

So  suppose that \( 27 a^2 + b^2 \) is nonzero and pick \( (\yy)\in \C^3 \) and \( (\zz)\in \C^3 \) on the line \(u_x= x_1+x_2+x_3=0 \). A convenient choice is  \( (\yy)=(1,\omega,\omega^2) \) and \( (\zz)=(1,\omega^2,\omega) \), because then \( f_{XXX} \), from \eqref{eqw97}, has a particularly simple form: \[ f_{XXX}=3\left((9 a - i \sqrt3\, b)X_1^3+(9 a +i \sqrt3\, b)X_2^3\right) \] 

If we further choose \( \alpha_1,\alpha_2\in\C \) such that \( \alpha_1^3=9 a - i \sqrt3\, b \) and \( \alpha_2^3=9 a +i \sqrt3\, b \), then \( f_{XXX} \) factors completely:
\[ f_{XXX}=3\left(\alpha_1 X_1+ \alpha_2 X_2\right)\left(\alpha_1 X_1+ \omega\alpha_2 X_2\right)\left(\alpha_1 X_1+ \omega^2\alpha_2 X_2\right) \]

Now choose  \( \gamma \in \C\) such that \( \gamma^3=9c \).  Since  \( \alpha_1^3\alpha_2^3=3(27a^2+b^2) \), the equation \( (27a^2+b^2)c+d^3=0 \) becomes  \( -\alpha_1^3\alpha_2^3 \gamma^3 =(3d)^3  \). So we can actually choose \( \gamma \) so that \( -\alpha_1\alpha_2\gamma= 3d \). 


To express the factors of \( f_{xxx} \) using \eqref{eq73}  we also need to calculate \( f_{xyy} \), \( f_{xyz}   \) and \( f_{xzz} \):
\begin{align*}
f_{xyy}&=3 (9 a - i \sqrt3\, b) (x_1 + \omega^2 x_2 + \omega x_3)=3\alpha_1^3 (x_1 + \omega^2 x_2 + \omega x_3) \\
f_{xyz}&=9 d (x_1 + x_2 + x_3)=-3 \alpha_1\alpha_2\gamma (x_1 + x_2 + x_3) \\
f_{xzz}&=3 (9 a + i \sqrt3\, b) (x_1 + \omega x_2 + \omega^2 x_3) =3\alpha_2^3 (x_1 + \omega x_2 + \omega^2 x_3)
\end{align*}

From  \eqref{eq73}, with \( A_1=\alpha_1 \), \( A_2=\alpha_2 \),  \( B_1=\alpha_1 \), \( B_2=\omega\alpha_2 \),  \( C_1=\alpha_1 \), and  \( C_2=\omega^2\alpha_2 \), the linear factors of \( f_{xxx} \) are (after removing the nonzero factor \( 3\alpha_1^2\alpha_2^2 \) from each),
\begin{equation}\label{eqq114}
\begin{aligned}
A_x&=\alpha_1(x_1 + \omega^2 x_2 + \omega x_3)+\gamma(x_1 + x_2 + x_3)+\alpha_2(x_1 + \omega x_2 + \omega ^2x_3) \\
B_x&=\alpha_1( \omega^2x_1 + \omega x_2 +  x_3)+\gamma(x_1 + x_2 + x_3)+\alpha_2(\omega x_1 + \omega ^2 x_2 + x_3) \\
C_x&=\alpha_1(\omega x_1 +  x_2 + \omega^2 x_3)+\gamma(x_1 + x_2 + x_3)+\alpha_2(\omega ^2x_1 +  x_2 + \omega x_3)
\end{aligned}
\end{equation}

Now one can check directly that \( 9f_{xxx}=A_xB_xC_x \) follows from the assumptions  \( \alpha_1^3=9 a - i \sqrt3\, b \),  \( \alpha_2^3=9 a +i \sqrt3\, b \), \( -\alpha_1\alpha_2\gamma= 3d \)  and  \( \gamma^3=9c \).
So, even though we found this factorization of \( f_{xxx} \) by assuming \( 27 a^2 + b^2 \) is nonzero, the validity of \( 9f_{xxx}=A_xB_xC_x \) does not depend on that assumption. 

\begin{theorem}
Suppose that \( f_{xxx} \) is given by \eqref{eqq113}. Then \( f_{xxx} \) is completely reducible if and only if \( a=b=c=0 \) or \( (27a^2+b^2)c+d^3=0 \). In the first case, the factorization of \( f_{xxx} \) follows  from \eqref{eqq117}. In the second case, the factorization of \( f_{xxx} \) follows from \( 9f_{xxx}=A_xB_xC_x \) where \( A_x \), \( B_x \), and \( C_x \) are given by \eqref{eqq114} and \( \alpha_1 , \alpha_2, \gamma \in \C\) satisfy
\[  \alpha_1^3=9 a - i \sqrt3\, b\qquad \alpha_2^3=9 a +i \sqrt3\, b\qquad  -\alpha_1\alpha_2\gamma= 3d \qquad\gamma^3=9c \] 
\end{theorem}

In retrospect we notice that \( A_x \), \( B_x \) and \( C_x \) are cycled among themselves by even permutations of \( \xx \). This means that \( f_{xxx} \) has the form \[ f_{xxx}= (r_ 1 x_ 1 + r_ 2 x_ 2 + r_ 3 x_ 3) (r_ 2 x_ 1 + r_ 3 x_ 2 + 
   r_ 1 x_ 3) (r_ 3 x_ 1 + r_ 1 x_ 2 + r_ 2 x_ 3) \] for some \( r_1,r_2,r_3\in \C \). It turns out that,  if \( c\neq 0 \), then \( r_1 \), \( r_2 \) and \( r_3 \) can be found directly as follows.
   
Supposing \( c\neq0 \), let \( r_1 \), \( r_2 \) and \( r_3 \) be the roots of the cubic polynomial  \[ F(x)= 27 c\, (x^3 - x^2) + (9 c + 3 d) x - (2 a + c + d).\]
The discriminant of \( F \) written in terms of its coefficients is \( -2^23^6 c (27 a^2 c + d^3)  \).  Using  \( d^3=-c(27a^2+b^2) \),  this can be written as \( 2^23^6 b^2 c^2 \). The discriminant of \( F \) written in terms of its roots is  \( 3^{12} c^4 (r_1 - r_2)^2 (r_1 - r_3)^2 (r_2 - r_3)^2 \). Equating these two expressions for the discriminant  and cancelling \( 3^6 c^2 \) gives 
\[ 2^2 b^2=3^6 c^2 (r_1 - r_2)^2 (r_1 - r_3)^2 (r_2 - r_3)^2\] 
Taking square roots of both sides of this equation, we have, after a possible reindexing of the roots of \( F \), 
\[ 2 b=27  c (r_1 - r_2) (r_1 - r_3)(r_2 - r_3)\] 
This equation, together with the equations expressing the coefficients of \( F \) in terms of the roots of \( F \), can be used to show that 
\[ f_{xxx}=27 c(r_ 1 x_ 1 + r_ 2 x_ 2 + r_ 3 x_ 3) (r_ 2 x_ 1 + r_ 3 x_ 2 + 
   r_ 1 x_ 3) (r_ 3 x_ 1 + r_ 1 x_ 2 + r_ 2 x_ 3) \]
More details on this derivation can be found in \cite{brookfield}.

Finally, we notice that the factors of \( f_{xxx}= x_1^3+x_2^3+x_3^3-3 x_1 x_2 x_3\) seen in \eqref{eqq117} are \textbf{not} cycled amongst themselves by even permutations of \( \xx \). This explains why this particular cubic form is a special case among cubic forms that are unchanged by even permutations of \( \xx \).

%

\end{document}